\documentclass{svmult}
\usepackage{amsmath}
\usepackage{amssymb,graphicx}
\usepackage[dvips]{epsfig}
\begin{document}

\title*{Stochastic Integrate and Fire Models: a
review on mathematical methods and their applications}
\author{Laura Sacerdote and Maria Teresa Giraudo}
\institute{Laura Sacerdote and Maria Teresa Giraudo \at Department of Mathematics, University of Torino, Via Carlo Alberto 10 Torino, Italy \email{laura.sacerdote@unito.it, mariateresa.giraudo@unito.it} \\ Work supported by MIUR PRIN 2008. LS is grateful to the organizers of the Summer School ''Stochastic Differential Equation Models with Applications to the Insulin-Glucose System and Neuronal Modelling'' for their kind hospitality.}

\titlerunning{Leaky Integrate and Fire models}
\maketitle

\textbf{\large Abstract}

Mathematical models are an important tool for neuroscientists. During the
last thirty years many papers have appeared on single neuron description and
specifically on stochastic Integrate and Fire models. Analytical results
have been proved and numerical and simulation methods have been developed
for their study. Reviews appeared recently collect the main features of
these models but do not focus on the methodologies employed to obtain them.
Aim of this paper is to fill this gap by upgrading old reviews on this
topic. The idea is to collect the existing methods and the available
analytical results for the most common one dimensional stochastic Integrate
and Fire models to make them available for studies on networks. An effort to
unify the mathematical notations is also made. This review is divided in
two parts:

1. Derivation of the models with the list of the available closed forms
expressions for their characterization;

2. Presentation of the existing mathematical and statistical methods for
the study of these models.

\section{Introduction}

\label{Section:Introd}

Progresses in experimental techniques, with the possibility to record
simultaneously from many neurons, move the interest of scientists from
single neuron to small or large networks models. Hence the time seems ripe
to summarize the contribution of single neuron models to our knowledge of
neuronal coding. Various types of spiking neuron models exist, with
different levels of details in the description. They range from biophysical
ones on the lines of the classic paper of 1952 by Hodgkin and Huxley \cite%
{HH}, to the "integrate and fire" variants (see for example \cite{gerstner}, 
\cite{diesmann}, \cite{plesser2}). Integrate and Fire (IF) type models
disregard biological details, that are accounted for through a stochastic
term, to focus on causal relationships in neuronal dynamics. Their relative
simplicity make them good candidates for the study of networks. Recent
reviews discuss qualitative (cf. \cite{Burkitt1}, \cite{Burkitt2}) and
quantitative (cf. \cite{BiolCybRev}) features of stochastic Leaky Integrate
and Fire (LIF) models. These models are a variant of IF
models where the spontaneous membrane decay is introduced. An older paper (%
\cite{Ricciardi}) concerns mathematical methods for their study. The aim of
this work is to collect the existing mathematical methods for LIF models, to
provide a set of methodologies for future studies on networks. Indeed,
although the stochastic LIF models are simplified representations of the
cells, they are considered good descriptors of the neuron spiking activity
(cf. for example \cite{jolivet}, \cite{kistler}). Though some criticisms
have appeared, showing some lacks in the fit of experimental data (cf. \cite%
{BiolCybRev}), these models are still largely employed. The most used is the
Ornstein-Uhlenbeck (OU) version but all of them have played a role for the understanding of the
mechanisms involved in neuronal code

The first IF models date back to 1907, when Lapique (\cite{Lapique})
proposed to describe the membrane potential evolution of a neuron, subject
to an input, using the time derivative of the law for the capacitance. In
the presence of an input current, the membrane voltage increases until it
reaches a constant threshold $S$. Then a spike occurs and the voltage is
reset to its resting value, to start again to evolve (cf. \cite{Tuckwell1}).
Although it reasonably fitted some experimental data this model lied
disregarded till the second half of the last century. Then it became the
embrional idea for "integrate and fire" models. The leading idea in the
formulation of stochastic IF and LIF models was to partition the features of
the neuron in two groups: the first ones were accounted for by the
mathematical description of the neuronal (deterministic) dynamics and the
second ones globally considered by means of a noise term. In Sect. \ref%
{Section:LIF} we derive the most popular LIF models after a brief
description of the biological features of interest in Sect. \ref%
{Section:Biol}.

Improvements of LIF models were proposed in the eighties, following the
initial illusion to become able to recognize the main laws governing our
brain. The lack of suitable mathematical instruments made soon clear the
difficulty to determine explicit expressions for the input-output
relationship. The end of the eighties and the starting years of the nineties
are characterized by mathematical and numerical advances, accompanied by the
development of new faster computers. Section \ref{Section:Math} is devoted
to a review of the main mathematical methods for the study of stochastic LIF
models, updating previous reviews \cite{Abrahams}, \cite{ricSato} and \cite%
{Ricciardi}.

In the nineties the use of such methodologies, as well as specific reliable
and powerful simulation methods, allowed to obtain a deeper knowledge of the
models' features. Unexpected results on the role of noise in neuronal coding
have been proved mathematically and confirmed experimentally (cf. for
example \cite{Segundo}).

Surprisingly, all research on LIF models has disregarded for a long time
their ability to fit real data. The only exception was \cite{LanskyStat}
that considered the parameter estimation problem. Recently papers on the
statistical estimation of model parameters started to appear. Section \ref%
{Section:Est} is dedicated to this subject.

\section{Biological features of the neuron}

\label{Section:Biol}

A comprehensive description of the physiological properties of neurons is
outside the aims of this work. We refer to \cite{Tuckwell1}, \cite{Tuckwell2}
and \cite{gerstner} for an exhaustive exposition of neurobiological
properties relevant in the modeling context.

The neurons are the elementary processing units in the central nervous
system, interconnected with intricate patterns. Neurons of different sizes
and shape, but sharing some fundamental features, exist in all the areas of
the brain. Their estimated number in the human brain is around $10^{12}$. A
typical neuron can be divided into three distinct parts called dendrites,
axon and soma. The dendrites play the role of the input device collecting
signals from other neurons and transmitting them to the soma. The soma is
the non-linear processing unit of the neuron. It generates a signal, known
as spike or action potential, if the total amount of inputs exceeds a
certain threshold. The axon is the output device carrying the signal to the
other neurons.

The action potentials are electrical pulses, having a duration of about $1-2$
$ms$ and an amplitude of around $100$ $mV$. They do not change their shape
along the transmission. A neuron cannot elicit a second spike immediately
after a first one has occurred due to the presence of a refractory period. A
chain of action potentials emitted by a single neuron is called a spike
train, a series of similar events occurring either at regularly spaced
instants of time or more randomly. The time between two consecutive spikes
is called interspike interval (ISI).

The site where the axon of a presynaptic neuron is linked with the dendrite
or the soma of a postsynaptic cell is the synapse, which often in the
vertebrates is of chemical type. When an action potential reaches a synapse,
it triggers complex bio-chemical reactions leading to the release of a
neurotransmitter and the opening of specific ionic channels on the membrane.
The ion influx leads to a change in the potential value at the postsynaptic
site and the translation of the chemical signal into an electrical one. This
voltage response is called the postsynaptic potential (PSP).

The effect of a spike on the postsynaptic neuron is measured in terms of the
potential difference between the interior of the cell and its surroundings,
called membrane potential. In the absence of spike inputs, the cell is at a
resting level of about $-65$ $mV$. If the change in membrane potential is
positive the synapse is excitatory and induces a negative depolarization,
otherwise it is inhibitory and hyperpolarizes the cell. In the absence of
inputs, i.e. in the silent state, the neuron membrane potential decays
exponentially toward the resting level.

The dimensions and number of synapses vary for different neurons. Some
neurons, such as Purkinje cerebellar cells, pyramidal neurons and
interneurons recorded in vitro (cf. for example \cite{BiolCybRev}), have a
huge number of synapses and extended dendritic trees. Integrate and fire
models can then be employed for the description of their output behavior
since, due to the large number of synapses, limit theorems can be used (cf. 
\cite{BiolCybRev}, \cite{Ricc}).

\bigskip

\section{One dimensional Stochastic Integrate and Fire Models}

\label{Section:LIF}

\subsection{Introduction and Notations}

\label{Subsection:Notations} The huge number of synapses impinging on the
neuron determines a stochasticity in the activating current not considered
in the Lapique model. The first attempt to formulate a stochastic IF model
is due to Gerstein and Mandelbrot. In \cite{gerstein} they fitted a
number of recorded ISI's through the Inverse Gaussian (IG) distribution, i.e. the
first passage time distribution of a Wiener process through a constant
boundary $S$. They described the membrane potential dynamics preceeding the
release of a spike through a Wiener process. To get a renewal process they
assumed that after each spike the membrane potential is instantaneously
reset to its initial value (cf. for example \cite{CoxMiller} for an
introduction on these processes). This model is the basis of successive more
realistic models.

In the models, classified as stochastic IF or LIF, one describes the time
evolution of the membrane potential by means of a suitable stochastic
process $X=\left\{ X_t ,t>t_0\right\}$ with $X_{t_0}=x_0$ and identifies the
ISI's with the random variable (r.v.) first passage time (FPT) of $X$ through the
threshold $S$:

\begin{equation}
T=T_S=\inf \left\{ t>t_0:X_t>S\right\}.  \label{FPT}
\end{equation}
The probability density function (pdf) of $T$, when it exists, is

\begin{equation}
g(t)=g(S,t\left|x_0,t_0\right.)=\frac{\partial}{\partial t} P(T<t).
\label{pdf}
\end{equation}
When $t_0=0$ we simply write $g(S,t\left|x_0\right.)$. In some instances $%
S=S(t)$.

In Subsections $3.2$, $3.5$ and $3.6$ we focus on models that describe the
subthreshold membrane potential as a diffusion process. In Subsects. $3.3$
and $3.4$ we present two continuos time Markov models, the Randomized Random Walk and Stein's model. Reviews on IF and LIF models have
already appeared (cf. \cite{ricSato}, \cite{Ricciardi}, \cite{Lanskysato})
but here we unify the notations and we list in a single contribution the
mathematical results sparse in different papers.

In the case of models using a diffusion process $X=\left\{ X_t
,t>t_0\right\} $, the diffusion interval is $I=(l,r)$, the drift coefficient
and infinitesimal variance (the infinitesimal moments) are:

\begin{eqnarray}
\mu(x)=\lim_{\Delta t\rightarrow 0} \frac{1}{\Delta t} E\left( \Delta
X_{t}\left\vert X_{t}=x\right. \right)  \label{def_infMoments} \\
\sigma^2(x)=\lim_{\Delta t\rightarrow 0} \frac{1}{\Delta t} E\left(
\left(\Delta X_{t}\right)^2\left\vert X_{t}=x\right. \right),  \notag
\end{eqnarray}
with $\Delta X_t=X_{t+\Delta t}-X_t$. The transition pdf $f\left( x,t\left|
x_{0},t_{0}\right. \right)=\frac{\partial P\left(X_t\leq x \left|X_{t_0}=x_0
\right.\right) }{\partial x}$ is solution of the Kolmogorov equation(cf. 
\cite{Ricc}):

\begin{equation}
\frac{\partial f\left( x,t\left\vert x_{0},t_{0}\right. \right) }{\partial
t_0}+\mu(x_0) \frac{\partial f\left( x,t\left\vert x_{0},t_{0}\right. \right)%
}{\partial x_0} +\frac{\sigma^2(x_0)}{2}\frac{\partial ^{2}f\left(
x,t\left\vert x_{0},t_{0}\right. \right) }{\partial x^{2}_0}=0  \label{Kolm}
\end{equation}
and of the Fokker-Planck equation 
\begin{equation}
\frac{\partial f\left( x,t\left\vert x_{0},t_{0}\right. \right) }{\partial t}%
=-\frac{\partial}{\partial x}\left\{\mu(x)f\left( x,t\left\vert
x_{0},t_{0}\right. \right)\right\} +\frac{1}{2}\frac{\partial ^{2}}{\partial
x^2}\left\{\sigma^2(x)f\left( x,t\left\vert x_{0},t_{0}\right.
\right)\right\}  \label{FP}
\end{equation}

with initial delta condition

\begin{equation}
\lim_{t_0\rightarrow t} f\left( x,t\left\vert x_{0},t_{0}\right.
\right)=\delta(x-x_0).  \label{d_wiener}
\end{equation}
Here $\delta$ denotes the Dirac delta function.

We suppose that the infinitesimal moments verify some mild conditions (cf. 
\cite{Karlin}, \cite{Ricc}, \cite{Ricciardi}) to guarantee the existence of
the solutions of the Fokker Planck and Kolmogorov equations. Furthermore,
when a dependence of the diffusion coefficients from $t$ is not specified,
the processes are time homogeneous, i.e. their properties are invariant with
respect to time shifts. When eq. (\ref{Kolm}) is solved in the presence of an absorbing boundary in $x=S$, a further
absorption condition must be imposed:

\begin{equation}
\lim_{x\rightarrow S} f^a\left( x,t\left\vert x_{0},t_{0}\right. \right)=0.
\label{abs}
\end{equation}
where $f^a\left( x,t\left\vert x_{0},t_{0}\right. \right)=\frac{\partial}{\partial x}P(X_t<x, T_S>t|X_s=y)$ is the
corresponding transition pdf. To get renewal processes, $X$ is always reset
to $x_0$ after each spike.

To characterize a diffusion model, one can also make use of the Ito-type
stochastic differential equation (SDE) verified by the process (cf. \textbf{%
SUSANNE}).

Jump diffusion models, allowing to distinguish the effect of neuronal inputs
according with their frequency and their size, are presented in Subsection %
\ref{Subsection:JD}.

The role of the threshold shape is illustrated in Subsect. \ref%
{Subsection:Shapes}, and the most recently introduced IF models are surveyed
in Subsection \ref{Subsection:Further}.

To switch from the description of the spike times of the neuron to the count
of the number of spikes up to a given time $t$, we introduce, in Subsection %
\ref{Subsection:Return}, the return processes.

\subsection{Wiener Process Model}

\label{Subsection:Wiener}

Gerstein and Mandelbrot (cf. \cite{gerstein}) described the time evolution
of the subthreshold membrane potential through a Wiener process $\overline{W}%
_{t}$ characterized by infinitesimal moments%
\begin{equation}
\mu \left( x\right) =\mu \text{ \ \ \ \ }\sigma ^{2}\left( x\right) =\sigma
^{2}  \label{diffusionWiener}
\end{equation}%
with $\mu \in \mathbb{R},:\sigma >0$. Their model was motivated by
experimental observations of the ISI's exhibiting histograms typical of
stable distributions. Indeed this property is exhibited by the FPT of a
Wiener process. One gets such
process from the standard Wiener process $W$ (cf. \textbf{SUSANNE})
through the transformation

\begin{equation}
\overline{W}_t =\mu t+\sigma W_t; \forall t\geq 0.  \label{Wstandard}
\end{equation}

To relate the use of the Wiener process with the membrane potential
evolution, Gerstein and Mandelbrot observed that the Wiener process is the
continuous limit of a random walk (cf. \textbf{SUSANNE}). The occurring of
jumps models the incoming of PSP's. The continuos limit is a good
approximation when the inputs are of small size and frequent . The
transition pdf of $W$ is

\begin{eqnarray}
f_{W}\left( x,t\left\vert x_{0},t_{0}\right. \right) &\equiv &\frac{\partial
P\left( W_{t}<x\left\vert W_{t_{0}}=x_{0}\right. \right) }{\partial x} \\
&=&\frac{1}{\sqrt{2\pi \sigma ^{2}\left( t-t_{0}\right) }}\exp \left\{ -%
\frac{\left[ x-x_{0}-\mu \left( t-t_{0}\right) \right] ^{2}}{2\sigma
^{2}\left( t-t_{0}\right) }\right\} .  \notag  \label{ddpWiener}
\end{eqnarray}%
To mimic the spiking times a constant absorbing boundary $S$ is introduced.
The spike times are then identified with the FPT, $T$, of the Wiener process
originated at $W_{t_{0}}=x_{0}$ through the boundary. To obtain the renewal
property, the process is instantaneously reset at $x_{0}$ after each spike.
Hence the ISI's correspond to the iid r.v.'s $T_{n}$, $n=1,2,...$,
with $T_{n}\sim T$.

The transition pdf of $W$, if $W_{t_0}=x_0$ is Gaussian with mean $%
E(W_t)=x_0+\mu t$ and variance $Var(W_t)=\sigma^2t$, while the FPT pdf
through a constant boundary $S>x_0$ is an IG
distribution, hence the pdf and the cumulative distribution are:%
\begin{equation}
g\left( S,t\left\vert x_{0}\right. \right) =\frac{S-x_{0} }{\sqrt{2\pi
\sigma ^{2}t^{3}}}\exp \left\{ -\frac{\left( S-x_{0}-\mu t\right) ^{2}}{%
2\sigma ^{2}t}\right\};  \label{FPTWiener}
\end{equation}

\begin{equation}
P(T<t)=\frac{1}{2}\left\{ Erfc\left[ \frac{S-x_{0}-\mu t}{\sigma \sqrt{2t}}%
\right] +e^{\frac{2\mu (S-x_{0})}{\sigma ^{2}}}Erfc\left[ \frac{S-x_{0}+\mu t%
}{\sigma \sqrt{2t}}\right] \right\} .  \label{FPTdistrWiener}
\end{equation}%
Here $Erfc$ denotes the complementary error function (cf. \cite{abramowitz}%
). The mean value and the variance of the FPT are

\begin{equation}
E(T)=\frac{S-x_0}{\mu}; \; \; Var\left(T\right)=\frac{(S-x_0) \sigma^2}{\mu^3%
}.  \label{ETwiener}
\end{equation}
The transition pdf in the presence of a constant absorbing boundary $S$ is
(cf. \cite{ricSato}):

\begin{eqnarray}
f^{a}(x,t|y,s)&=&\frac{1}{\sigma\sqrt{2\pi (t-s)}} \left(\exp[-\frac{(x-y-\mu (t-s))^2}{%
2\sigma^2(t-s)}]\right.  \notag \\
&-& \left.\exp[\frac{2\mu}{\sigma^2}(S-y)-\frac{(x-2S+y-\mu (t-s))^2}{%
2\sigma^2(t-s)}]\right).  \label{faTwiener}
\end{eqnarray}

Despite the excellent fitting with some experimental data, Gerstein and
Mandelbrot model was criticized for its biological simplifications (cf. \cite%
{Tuckwell2}). However it allows to obtain results that help the intuition for more
realistic models and it is still used for this aim, taking advantage of the
existence of a closed form FPT pdf through a constant boundary.

Its FPT pdf is known also through particular time dependent boundaries.
These FPT's can be used to account for the refractory period following a
spike. Indeed a time varying boundary, assuming high values at small times
and then decreasing, makes short ISI's rare. The FPT pdf is known for a
continuos piecewise-linear boundary $S(t)=\alpha _{i}+\beta _{i}t,:t\in
\lbrack t_{i-1},t_{i}],:i\geq 1$ where $t_{0}<t_{1}<t_{2}<...$ and $\alpha
_{i}$, $\beta _{i}\in R$ with $t_{0}\geq 0$. If $t\in \lbrack 0,\infty )$,
such boundary is linear and the FPT pdf is

\begin{equation}
g(\alpha _{1}+\beta _{1}t,t\left\vert x_{0}\right. )=\frac{\left\vert \alpha
_{1}-x_{0}\right\vert }{\sigma \sqrt{2\pi t^{3}}}\exp {-\frac{\left[ \alpha
_{1}+\beta _{1}t-\mu t-x_{0}\right] ^{2}}{2\sigma ^{2}t}}.  \label{linWiener}
\end{equation}%
In the general case, setting $\alpha _{i+1}=\alpha _{i}+\beta _{i}t_{i}$ one
gets that $t\mapsto S(t)$ is continuous on $[t_{0},\infty )$. If we put $%
S_{i}=S(t_{i})$, the transition pdf for the process $W$ in the presence of
absorbing boundary $S\left( t\right)$, $%
f^{a}(x_{1},t_{1};x_{2},t_{2};...;x,t|x_{0},t_{0})$, is for $t\in
(t_{n-1},t_{n})$ (cf. \cite{WangPotz}):

\begin{eqnarray}
f^{a}(x_{1},t_{1};x_{2},t_{2};...;x,t|x_{0},t_{0})= \frac{\partial
^{n}}  {\partial x_{1}...\partial x_{n}} \\ \notag
\left\{P\left( W_{t_{1}}<x_{1},...,W_{t_{n-1}}<x_{n-1},W_{t}<x; T>t \left\vert W_{t_0}=x_0<\alpha_1\right.\right)\right\} \\ \notag
\times \prod_{i=1}^{n-1}f^{a}(x_{i},t_{i}\left\vert x_{i-1},t_{i-1})f^{a}(x,t|x_{n-1},t_{n-1}\right.) \\ \notag
\times \prod_{i=1}^{n-1}\left( 1-e^{-2\frac{(S_{i}-x_{i})(S_{i-1}-x_{i-1})%
}{t_{i}-t_{i-1}}}\right) \\ \notag
\times \frac{f^{a}(x,t|x_{n-1},t_{n-1})}{\sqrt{2\pi (t_{i}-t_{i-1})}}\exp
\left( -\frac{(x_{i}-x_{i-1})^{2}}{2(t_{i}-t_{i-1})}\right) 
\label{WangPotz}
\end{eqnarray}
for $x_{i}\leq S_{i}$, $1\leq i\leq n$; $x_{0}<S_{0}$ and $x\in (-\infty ,S)$.

Further closed form expressions for the FPT of a Wiener process have been
obtained by the method of images (cf. \cite{Daniels2}) or as solutions of
suitable integral equations (\cite{RicSacSato}), when $\left\vert \frac{%
dS\left( t\right) }{dt}\right\vert \leq Ct^{-\alpha }$, with $\alpha <1/2$%
and $C$ a constant. This last case is discussed in Section \ref%
{Subsubsection:Int}. \ and involves series of multiple integrals.

\subsection{Randomized Random Walk Model}

\label{Subsection:RRW}

In the Randomized Random Walk (RRW) the regularly spaced intertimes of the
random walk between PSPs are substituted with exponentially distributed
intertimes of parameters $\lambda ^{+\text{ }}$and $\lambda ^{-}$ for
excitatory and inhibitory PSP's respectively. The process $X$ with $X_0=0$
has mean and variance

\begin{equation}
E\left( X_{t}\right) =\delta \left( \lambda ^{+}-\lambda ^{-}\right)
t;\;\;Var\left( X_{t}\right) = \delta ^{2}\left( \lambda ^{+}+\lambda
^{-}\right) t,  \label{RRW_moments}
\end{equation}
respectively. Here $\delta>0$ is the constant amplitude of PSP's. The FPT
pdf through the boundary $S$, with $S$ an integer multiple of $\delta$, is
(cf. \cite{Tuckwell2}): 
\begin{equation}
g\left( S,t\left\vert 0\right. \right) = \frac{S}{\delta} \left( \frac{%
\lambda ^{+}}{\lambda ^{-}}\right) ^{ S/2 \delta}\frac{e^{-\left( \lambda
^{+}+\lambda ^{-}\right) t}}{t}I_{S/ \delta }\left( 2t\sqrt{\lambda
^{+}\lambda ^{-}}\right), \text{ \ \ }t>0  \label{FPT_dist_RRW}
\end{equation}
where $I_{\eta}\left( .\right) $ is the modified Bessel function of
parameter $\eta$ (cf. \cite{abramowitz}). The mean and variance of the ISI
distribution are (cf. \cite{Tuckwell2}):

\begin{equation}
E(T) =\frac{S }{\delta \left( \lambda ^{+}-\lambda ^{-}\right) };\text{ \ \ }%
Var\left( T\right) =\frac{ S \left( \lambda ^{+}+\lambda ^{-}\right) }{%
\delta \left( \lambda ^{+}-\lambda ^{-}\right) ^{3}}.
\label{MeanVarianceRRW}
\end{equation}
When $\delta \rightarrow 0$, assuming $\lambda^+ \sim \lambda^- \sim \frac{1%
}{2 \delta^2}$, this model converges to the Wiener process with $\mu=0$.

\subsection{Stein's model}

\label{Subsection:Stein}

In 1965 Stein (cf. \cite{stein}) formulated the first LIF model, i.e. an IF
model with the leakage feature, by introducing the spontaneous
membrane decay in the absence of PSP's in the RRW model. The process $%
X $ is solution of the SDE 
\begin{equation}
dX_{t}=\left( -\frac{X_{t}}{\theta }+\rho \right) dt+\delta
^{+}dN_{t}^{+}+\delta ^{-}dN_{t}^{+};\;X_{t_{0}}=x_{0}.  \label{Stein_eq}
\end{equation}%
Here $\theta >0$ is the membrane time constant, $\rho $ is the resting
potential, $\ N_{t}^{+}$ and $N_{t}^{-}$ are independent Poisson processes
of parameters $\lambda ^{+}$ and $\lambda ^{-}$ respectively and $\delta
^{+}>0$, $\delta ^{-}<0$ are the amplitudes of excitatory and inhibitory
PSP's. Generally for this model and for all those descending from it one
assumes $\rho =0,$ since the case $\rho \neq 0$ can be obtained by the
shift $X\mapsto X-\rho $. Following the IF models structure, the
spike times are the first crossing times of the process
through the boundary and the membrane potential is instantaneously reset to
its resting value after each spike. The infinitesimal moments are:%
\begin{eqnarray}
M_{1}\left( x\right) &=&\lim_{h\rightarrow 0}\frac{E\left[
X_{t+h}-X_{t}\left\vert X_{t}=x\right. \right] }{h}=-\frac{x}{\theta }+\rho
+\lambda ^{+}\delta ^{+}+\lambda ^{-}\delta ^{-}  \label{momStein} \\
M_{2}\left( x\right) &=&\lim_{h\rightarrow 0}\frac{E\left[ \left(
X_{t+h}-X_{t}\right) ^{2}\left\vert X_{t}=x\right. \right] }{h}=\lambda
^{+}\left( \delta ^{+}\right) ^{2}+\lambda ^{-}\left( \delta ^{-}\right)
^{2}.  \notag
\end{eqnarray}%
The mean trajectory, in the absence of a threshold, is%
\begin{equation}
E\left( X_{t}\left\vert X_{0}=x_{0}\right. \right) =x_{0}e^{-t/\theta
}+\left( \lambda ^{+}\delta ^{+}+\lambda ^{-}\delta ^{-}\right) \theta
\left( 1-e^{-t/\theta }\right) .  \label{meanStein}
\end{equation}%
The FPT problem for the process (\ref{Stein_eq}) is still unsolved and the
use of simulation techniques is required for its analysis.

\subsection{Ornstein-Uhlenbeck Diffusion Model}

\label{Subsection:OU}

The OU process was proposed as a continuos limit of the
Stein model to facilitate the solution of the FPT problem. The rationale for
this limit is the huge number of synapses characterizing certain neurons
such as the Purkinje cells. The PSP's determine frequent small jumps and
limit theorems can be applied to get a diffusion process. The OU process,
already known in the Physics literature (cf. \cite{OrnUhl}), belongs to both
classes of Markov and Gaussian stochastic processes.

Different approaches can be followed to obtain the diffusion limit of Stein's
model. In \cite{Kallianpur} the convergence of the measure of a Stein's
process to that of the OU one is studied.In \cite{capocelli}, \cite%
{Ricciardi} it is proved that the continuos limit of the transition pdf for
the process (\ref{Stein_eq}) converges to a pdf that verifies the Fokker
Plank equation for the OU process. Alternatively the OU model can be derived
from a differential equation describing the membrane potential dynamics. We
sketch in the following these last two approaches.

Due to the time continuity and to the Markov property of Stein's process,
the Smolukowsky equation holds for the transition pdf:

\begin{equation}
f\left( x,t+\Delta t\left\vert x_{0},t_{0}\right. \right) =\int_{-\infty
}^{\infty }f\left( x,t+\Delta t\left\vert z,t\right. \right) f\left(
z,t\left\vert x_{0},t_{0}\right. \right) dz.  \label{Smolukowski}
\end{equation}

In the absence of inputs to the neuron, the process (\ref{Stein_eq}),
initially in the state $z$ at time $t$, decays exponentially to the zero
resting potential, reaching at time $t+\Delta t$ the value $%
x_{1}=ze^{-\Delta t/\theta }$. In the case of an excitatory or an inhibitory
input at time $u\in \left( t,t+\Delta t\right) $ the potential becomes $%
x_{2}\left( u\right) =ze^{-\Delta t/\theta }+\delta ^{+}e^{-\left( t+\Delta
t-u\right) /\theta }$ or $x_{3}\left( u\right) =ze^{-\Delta t/\theta
}+\delta ^{-}e^{-\left( t+\Delta t-u\right) /\theta }$, respectively.
Setting $x_{i}=\frac{1}{\Delta t}\int_{t}^{t+\Delta t}x_{i}\left( u\right)
du,i=2,3$ the l.h.s. of (\ref{Smolukowski}) becomes:

\begin{eqnarray}
f\left( x,t+\Delta t\left\vert z,t\right. \right) &=&\left[ 1-\left( \lambda
^{+}+\lambda ^{-}\right) \Delta t\right] \delta \left( x-x_{1}\right)  \notag
\\
&+&\lambda ^{+}\Delta t\delta \left( x-x_{2}\right) +\lambda ^{-}\Delta
t\delta \left( x-x_{3}\right) +o(\Delta t).  \label{transition}
\end{eqnarray}%
Hence (\ref{Smolukowski}) becomes:

\begin{eqnarray}
f\left( x,t+\Delta t\left| x_{0},t_{0}\right. \right)&=&e^{\Delta t/\theta }
\left\{ \left[ 1-\left( \lambda ^{+}+\lambda ^{-}\right) \Delta t\right]
f\left( xe^{\Delta t/\theta },t\left| x_{0},t_{0}\right. \right)\right. 
\notag \\
&+&\left. \lambda ^{+}\Delta t f \left(xe^{\Delta t/\theta }-\frac{\delta
^{+}\theta }{\Delta t}\left( e^{\Delta t/\theta }-1\right) ,t\left|
x_{0},t_{0}\right. \right) \right.   \\
&+&\left. \lambda ^{-}\Delta t f(xe^{\frac{\Delta}{\theta}}-\frac{\delta
^{-}\theta }{\Delta t}( e^{\frac{\Delta}{\theta}}-1),t\left| x_{0},t_{0}
\right.)\right\}+o( \Delta t).  \notag
\label{develop_f}
\end{eqnarray}
Approximating $e^{\Delta t/\theta }\approx \frac{\Delta t}{\theta}+1$,
dividing by $\Delta t$, when $\Delta t\rightarrow 0,$ (\ref{develop_f})
becomes

\begin{eqnarray}
\frac{\partial f\left( x,t\left\vert x_{0},t_{0}\right. \right) }{\partial t}
=\frac{\partial }{\partial x}\left( \frac{x}{\theta }f\left( x,t\left\vert
x_{0},t_{0}\right. \right) \right) + \lambda ^{+} \left[ f\left( x-\delta
^{+},t\left\vert x_{0},t_{0}\right. \right)\right.  \notag \\
\left.-f\left( x,t\left\vert x_{0},t_{0}\right. \right) \right] + \lambda
^{-}\left[ f\left( x-\delta ^{-},t\left\vert x_{0},t_{0}\right. \right)
-f\left( x,t\left\vert x_{0},t_{0}\right. \right) \right].  \label{eq_f}
\end{eqnarray}
Developing the terms in square brackets as Taylor series around $x$

\begin{eqnarray}
\frac{\partial f\left( x,t\left\vert x_{0},t_{0}\right. \right) }{\partial t}
&=&-\frac{\partial }{\partial x}\left[ -\left( \frac{x}{\theta }+\lambda
^{+}\delta ^{+}+\lambda ^{-}\delta ^{-}\right) f\left( x,t\left\vert
x_{0},t_{0}\right. \right) \right] +  \label{q_f2} \\
&&\sum_{n=2}^{\infty }\frac{\left( -1\right) ^{n}}{n!}\frac{\partial ^{n}}{%
\partial x^{n}}\left\{ \left[ \left( \delta ^{+}\right) ^{n}\lambda
^{+}+\left( \delta ^{-}\right) ^{n}\lambda ^{-}\right] f\left( x,t\left\vert
x_{0},t_{0}\right. \right) \right\} .  \notag
\end{eqnarray}%
Assuming $\delta ^{+}=-\delta ^{-}=\delta $, $\lambda ^{+}\sim \frac{A^{+}}{%
\delta }+\frac{\sigma ^{2}}{2\delta ^{2}}$, $\lambda ^{-}\sim \frac{A^{-}}{%
\delta }+\frac{\sigma ^{2}}{2\delta ^{2}},$ with $\sigma ,A^{+},A^{-}$
positive constants, for $\delta \rightarrow 0$ we get:

\begin{equation}
\frac{\partial f\left( x,t\left\vert x_{0},t_{0}\right. \right) }{\partial t}%
=\frac{\partial }{\partial x}\left[ \left( -\frac{x}{\theta }+\mu \right)
f\left( x,t\left\vert x_{0},t_{0}\right. \right) \right] +\frac{\sigma ^{2}}{%
2}\frac{\partial ^{2}f\left( x,t\left\vert x_{0},t_{0}\right. \right) }{%
\partial x^{2}},  \label{F_P_eq}
\end{equation}
i.e. the Fokker Plank equation for an OU process with $\mu =A^{+}-A^{-}.$
Its solution with an initial delta condition (\ref{d_wiener}) is the
transition pdf

\begin{eqnarray}
f_{OU}\left( x,t\left\vert x_{0},t_{0}\right. \right) &\equiv& \frac{%
\partial P\left( X_{t}<x\left\vert X_{t_{0}}=x_{0}\right. \right) }{\partial
x} = \frac{1}{\sigma\sqrt{\pi \theta \left( 1-e^{-\frac{2\left(
t-t_{0}\right)}{\theta} }\right) }}  \notag \\
&\times& \exp \left\{ -\frac{\left[ x-x_{0}e^{-\frac{\left( t-t_{0}\right)}{%
\theta} }-\mu \theta (1-e^{-\frac{\left( t-t_{0}\right)}{\theta} })\right]
^{2}}{\sigma ^{2}\theta \left( 1-e^{-\frac{2\left( t-t_{0}\right)}{\theta}
}\right) }\right\}.  \label{ddpOu}
\end{eqnarray}

The diffusion interval coincides with the real line and the mean and
variance of the OU process $X_{t}$ with $X_{t_0}=x_0$ are:

\begin{equation}
E\left( X_{t}|x_0\right) =\mu \theta \left(1-e^{-t/\theta }\right)+x_0
e^{-t/\theta } ;\; Var\left( X_{t}|x_0\right) =\frac{\sigma ^{2}\theta }{2}%
\left( 1-e^{-2t/\theta }\right).  \label{mean_varianceOU}
\end{equation}

Properties of the models, as well as the range of validity of some
approximate formulae for the FPT problem, depend upon the value of the
asymptotic mean depolarization of the process $X$. Hence in various
instances we distinguish between two distinct firing regimes, subthreshold
if $E(X_{\infty })<S$ and suprathreshold in the opposite case.

\begin{figure}[htp]
\centering
\includegraphics[width=11cm]{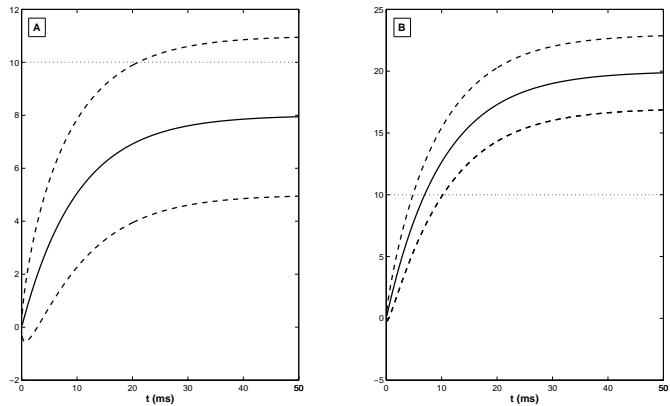}
\caption{Mean value (\textit{middle line}) and curves $E(X_t)-3Var(X_t)$ (%
\textit{lower line}) and $E(X_t)+3Var(X_t)$ (\textit{upper line}) for an OU
process with parameters $\protect\mu=0.8$ $mVms^{-1}$, $\protect%
\sigma^2=0.2$ $mV^2ms^{-1}$ and $\protect\theta=10$ $ms$.}
\label{fig:MeanVar}
\end{figure}

Model of LIF type can be interpreted in the frame of threshold detectors
theory. The presence of a feeble noise helps the detection of the signal, a
characteristic of any threshold detector. In Fig. \ref{fig:MeanVar} we plot
the mean value of an OU process (\ref{mean_varianceOU}a) together with $%
E(X_{t})-3Var(X_{t})$ and $E(X_{t})+3Var(X_{t})$, making use of (\ref%
{mean_varianceOU}b). The two panels correspond to examples in the
subthreshold regime (Panel A) and in the suprathreshold regime (Panel B).
The intrinsic random variability determines crossings even in the
subthreshold regime.

The OU model can also be obtained from the differential equation for the
time evolution of the subthreshold membrane potential in the presence of
spontaneous decay of parameter $\theta$ and net input $\mu$:

\begin{equation}
\frac{dX_{t}}{dt} =-\frac{X_{t}}{\theta }+\mu ; \; X_0 =x_{0}.
\label{differential_leaky}
\end{equation}
Adding a noise term of intensity $\sigma$ to account for the random PSP's,
one gets:

\begin{equation}
dX_{t}=\left( -\frac{X_{t}}{\theta }+\mu \right) dt+\sigma
dW_{t};\;X_{0}=x_{0}.  \label{SDE_OU}
\end{equation}%
This is the SDE of an OU process (cf. \textbf{SUSANNE}). The analytical
expression for the FPT pdf of the OU process is still an open problem. In 
\cite{Alili} three alternative representations of the distribution of $T$
are introduced. The first one involves an eigenvalue expansion in terms of
zeros of the parabolic cylinder functions while the second is an integral
representation involving special functions. We report here the third one
that writes the FPT pdf of an OU process with $\mu =0$ through the boundary $%
S$ in terms of a three-dimensional Bessel bridge:

\begin{equation}
g(t)=e^{-(S^{2}-x_{0}^{2}-t)/2\theta }g_{W}(t)E_{Bb}\left[ \exp \left( -%
\frac{1}{2\theta ^{2}}\int_{0}^{t}(r_{s}-S)^{2}ds\right) \right] .
\label{FPTOUBB}
\end{equation}%
Here $g_{W}(t)$ is the FPT pdf through the boundary $S$ for the standard
Wiener process, $r_{s}$ is the three-dimensional Bessel bridge over the
interval $[0,t]$ between $r_{0}=0$ and $r_{t}=S-x_{0}$. This process is
solution of:

\begin{equation}
dr_{s}=\left( \frac{y-r_{s}}{t-s}+\frac{1}{r_{s}}\right)
ds+dW_{s},\;r_{0}=x,s<t.  \label{BesselBridge}
\end{equation}%
In (\ref{FPTOUBB}) $E_{Bb}$ indicates the expectation with respect to the
Bessel bridge law. In \cite{ZhangFeng} formula (\ref{FPTOUBB}) is used to
approximate the FPT pdf with Monte Carlo techniques.

An explicit expression for the FPT density of continuous Gaussian processes
to a general boundary is obtained under mild conditions in \cite{durbin2},
while an expression for the FPT of the Wiener process to a curved boundary
is expanded as a series of multiple integrals in \cite{durbin3}.

Existing available closed form expression include the case of a hyperbolic
boundary (cf. \cite{buonocore})

\begin{equation}
S(t)=\mu \theta +Ae^{-\frac{t}{\theta }}+Be^{\frac{t}{\theta }},
\label{hyperbolic}
\end{equation}%
with $A$ and $B$ arbitrary constants (cf. \cite{Ricc}). Furthermore specific
boundaries can be obtained through the space time transformations described
in Section \ref{Subsubsection:Change} applied to closed form solutions for
the case of the Wiener process. The Laplace transform of the FPT pdf in the
case of a constant boundary $S$ is (cf. \cite{RicSac}) :

\begin{equation}
E\left( e^{-\lambda T}\right) =\exp \left\{ \frac{\left( x_{0}-\mu \theta
\right) ^{2}-\left( S-\mu \theta \right) ^{2}}{2\sigma ^{2}\theta }\right\} 
\frac{D_{-\lambda \theta }\left[ \sqrt{\frac{2}{\sigma ^{2}\theta }}\left(
\mu \theta -x_{0}\right) \right] }{D_{-\lambda \theta }\left[ \sqrt{\frac{2}{%
\sigma ^{2}\theta }}\left( \mu \theta -S\right) \right] },  \label{Laplace}
\end{equation}%
where $D_{\nu }\left( .\right) $ is the Parabolic Cylinder Function (cf. 
\cite{abramowitz}) of parameter $\nu $. No analytical inversion formula is
available for eq. (\ref{Laplace}). Reliable and efficient procedures,
discussed in Sect. \ref{Section:Math}, can be applied to obtain the FPT pdf
either numerically or by means of simulations for constant or time dependent
boundaries.

The FPT mean has been determined as derivative of (\ref{Laplace}), computed
in $\lambda =0$ (cf. \cite{RicSac}):

\begin{eqnarray}
E(T) &=&\theta \left\{\frac{1}{2} \left( \sum_{n=1}^{\infty }\frac{{\small x}%
_{S}^{2n}}{{\small n}\left( 2n-1\right) {\small !!}} -\sum_{n=1}^{\infty }%
\frac{x_{1}^{2n}}{{\small n}\left( 2n-1\right) {\small !!}} \right) \right. 
\notag \\
&&\left. +\sqrt{\pi /2}\left[x_{1}\phi \left( \frac{1}{2}{\small ,}\frac{3}{2%
}{\small ;}\frac{x_{1}^{2}}{2}\right) {\small -}x_{S}\phi \left( \frac{1}{2}%
\frac{3}{2}\frac{x_{S}^{2}}{2}\right) \right] \right\}  \label{ETou}
\end{eqnarray}%
where $x_{1}=(\mu \theta - x_{0})\sqrt{2/\left( \sigma ^{2}\theta \right) }$%
, $x_{S}=(\mu \theta - S)\sqrt{2/\left( \sigma ^{2}\theta \right)}$ and $%
\phi \left( a,c;z\right)$ is the Kummer function (cf. \cite{abramowitz}).
Alternatively the mean is expressed through the Siegert formula (cf. \cite%
{siegert}): 
\begin{equation}
E(T)=\sqrt{\frac{\pi \theta}{\sigma^2}}\int^{S-\mu \theta}_{-\mu
\theta}\left\{1+Erf(\frac{z}{\sigma \sqrt{\theta}})\right\}\exp(\frac{z^2}{%
\sigma^2 \theta})dz,  \label{momOUSiegert}
\end{equation}
where $Erf(.)$ denotes the error function (cf. \cite{abramowitz}).

Use of (\ref{ETou}) or (\ref{momOUSiegert}) depends on the value of the
parameters since the two formulae present numerical difficulties for
different ranges. Approximate formulae (cf. \cite{LanskySac}) hold for
specific ranges. If $\mu \theta >S$ and $\sigma \rightarrow 0$, i.e. in the
quasi-deterministic case, the mean FPT\ can be approximated by equating the
expression of $E(X_t|x_{0})$ with $S$ to obtain (cf. \cite{tomasSac2}):

\begin{equation}
E(T)\approx -\theta \ln \left(\frac{S-\mu \theta}{x_{0}-\mu \theta}\right).
\label{ETouap1_0}
\end{equation}%
Note that (\ref{ETouap1_0}) disregards the effect of the noise on the
crossings. If $x_{0}<<S,$ or equivalently if $\sigma $ is sufficiently small
and $\mu $ is negative so that the crossing is rare event, the
approximation

\begin{equation}
E(T)\approx \frac{\sigma \sqrt{\pi \theta ^{3}}}{S-\mu \theta }\exp \left( 
\frac{(S-\mu \theta )^{2}}{\sigma ^{2}\theta }\right)  \label{ETouap2}
\end{equation}%
holds (cf. \cite{gioNoRiasymp}). A linear approximation for the firing rate $%
f=1/E(T)$, obtained using (\ref{ETou}), is (cf. \cite{LanskySac}):

\begin{equation}
f\left( \mu \right) =\frac{1}{\pi \theta S}\left( \sigma \sqrt{\pi \theta }%
+2\theta \mu -S\right).  \label{freq}
\end{equation}
This approximation holds when $\left( \mu \sqrt{\theta }\right) /\sigma$ and 
$\left( \mu \sqrt{\theta }-S \theta\right) /\sigma$ are small enough.

When neither (\ref{ETou}) nor (\ref{momOUSiegert}) are suitable for
computations and $\mu \tau \gg S$ but $\sigma$ is not small enough to apply
approximation (\ref{ETouap1_0}) an ''ad hoc'' procedure to evaluate the
mean FPT is possible. One establishes at first the time $t_1$ at which $%
E(X_t)+2\sqrt{Var(X_t)}$ crosses the threshold $S$, i.e. when most
trajectories are still below the threshold. For $t>t_1$ the process is then
approximated by means of the Wiener process with drift $\mu$ and initial
value $E(X_{t_1})$.

The second moment of the FPT for the OU process is (cf. \cite{nobricSac2}):

\begin{eqnarray}  \label{ET2ou}
E(T^2) &=& 2 \theta E(T) \left[\sqrt{\pi} \varphi_1(\frac{x_s}{\sqrt{2}}) +
\psi_1(\frac{x_s}{\sqrt{2}})\right] \\
&+& 2 \theta^2 \left\{\sqrt{\pi} \ln 2 \left[\varphi_1(\frac{x_s}{\sqrt{2}})
- \varphi_1(\frac{x_1}{\sqrt{2}})\right]\right.  \notag \\
&-& \left.\sqrt{\pi} \left[\varphi_2(\frac{x_s}{\sqrt{2}}) + \varphi_2(\frac{%
x_1}{\sqrt{2}})\right]-\psi_2(\frac{x_s}{\sqrt{2}})+\psi_2(\frac{x_1}{\sqrt{2%
}}) \right\}  \notag
\end{eqnarray}
where $x_{1}$, $x_{S}$ are defined as in (\ref{ETou}) and

\begin{eqnarray}
\varphi_1(z) &=& \int_0^z e^{t^2}dt=\sum^{\infty}_{k=0}\frac{z^{2k+1}}{%
k!(2k+1)};  \notag \\
\varphi_2(z) &=& \sum_{n=0}^{\infty}\frac{z^{2n+3}}{(n+1)!(2n+3)}%
\sum_{k=0}^{n}\frac{1}{2k+1}  \notag \\
\psi_1(z) &=& 2\int_0^z e^{u^2}\int_0^u e^{-v^2}dv du= \sum^{\infty}_{k=0}%
\frac{2k z^{2k+2}}{(2k+1)!!(k+1)};  \notag \\
\psi_2(z) &=& \sum_{n=0}^{\infty}\frac{2^n z^{2n+4}}{(2n+3)!!(n+2)}%
\sum_{k=0}^{n}\frac{1}{k+1}.
\end{eqnarray}

In \cite{cerbone} the mean, variance and skewness of the FPT for the OU
process are tabulated for neurobiologically compatible choices of the
parameters.

Asymptotic results for the FPT of the OU process are presented in Secion \ref%
{Subsubsection:Asym}.

\subsection{Reversal Potential Models}

\label{Subsection:RevPot}

The diffusion interval of the OU process is the real line but
large negative values of the membrane potential are unrealistic. Hence other
models introduce a saturation effect on the membrane sensibility. When the
value of the membrane potential is close to the reversal potential $V_{I}$
the incoming inputs produce a reduced effect (cf. \cite{Lanska}, \cite%
{giorno}). A diffusion model with reversal potentials is proposed in \cite%
{giorno} as a diffusion limit on a birth and death process. A similar
diffusion limit is obtained in \cite{Lanska} from a variant of Stein's model
(\ref{Stein_eq}) where an inhibitory reversal potential $V_{I}$ is
introduced:

\begin{equation}
dY_{t}=-\frac{1}{\theta }Y_{t}dt+\delta ^{+}dN_{t}^{+}+\left[ \varepsilon
\left( Y_{t}-V_{I}\right) +\xi \sqrt{Y_{t}-V_{I}}\right] dN_{t}^{-}; \;
Y_{0}=y_{0}.  \label{steinReversal}
\end{equation}
Here $N_{t}^{+}$, $N_{t}^{-}$ $\delta ^{+}$ and $\theta $ are the same as in
(\ref{Stein_eq}), $\varepsilon \in \left( -1,0\right) $, $V_{I}<x_{0\text{ }%
} $ are two constants and $\xi$ is a suitably defined random variable. The
first two infinitesimal moments (\ref{def_infMoments}) of this model are:

\begin{eqnarray}
\mu (y) &=&-y\left( \frac{1}{\theta }-\varepsilon \lambda ^{-}\right)
+\delta ^{+}\lambda ^{+}-\varepsilon \lambda ^{-}V_{I};  \label{infMoments} \\
\sigma ^{2}(y) &=&\lambda ^{+}\left( \delta ^{+}\right) ^{2}+\lambda
^{-}\varepsilon ^{2}\left( y-V_{I}\right) ^{2}+\lambda ^{-}Var\left( \xi
\right) \left( y-V_{I}\right) .  \notag
\end{eqnarray}%
The mean trajectory of the process originated in $Y_{0}=y_{0}$ is 
\begin{equation}
E\left( Y_{t}\left\vert y_{0}\right. \right) =x_{0}e^{-t\left( \frac{1}{%
\theta }-\varepsilon \lambda ^{-}\right) }+\frac{\lambda ^{+}\delta
^{+}-\varepsilon \lambda ^{-}V_{I}}{\frac{1}{\theta }-\varepsilon \lambda
^{-}}\left( 1-e^{-t\left( \frac{1}{\theta }-\varepsilon \lambda ^{-}\right)
}\right) .  \label{meanRev}
\end{equation}

The diffusion limit of this model (known in the neurobiological literature
as the Feller model and in other contexts as the Cox-Ingersol-Ross process)
is identified with the solution to the SDE (cf. \cite{Lanska})

\begin{equation}
dY_{t}=\left( -\frac{Y_{t}}{\tau }+\mu _{2}\right) dt+\sigma _{2}\sqrt{%
Y_{t}-V_{I}}dW_{t};\;Y_{0}=y_{0}.  \label{feller}
\end{equation}%
Here the constants $\mu _{2},\sigma _{2}$ and $\tau $ are related with those
of model (\ref{steinReversal}) by imposing the equality of the infinitesimal
moments. One first sets $\tau =\frac{\theta }{1-\varepsilon \lambda
^{-}\theta }$ and $\mu _{2}=\lambda ^{+}\delta ^{+}-\varepsilon \lambda
^{-}V_{I}=\lambda ^{+}\delta ^{+}-\lambda ^{-}\Delta _{I}^{-}$ where $\Delta
_{I}^{-}=\varepsilon V_{I}$, then substitutes the variable $\xi $ in (%
\ref{steinReversal}) by a suitable sequence of r.v.'s such that in
the limit for $n\rightarrow \infty $ one gets $Var(\xi _{n})=0$. This choice
allows to obtain the same infinitesimal variance at the resting level for
the two processes. Hence one gets

\begin{equation}
\sigma _{2}^{2}=-\frac{\lambda ^{+}(\delta ^{+})^{2}+\lambda ^{-}(\Delta
_{I}^{-})^{2}}{V_{I}}.  \label{limfel3}
\end{equation}%
Note that, due to the expressions for $\tau $ and $\mu _{2}$, the parameters 
$\tau $ and $\mu _{2}$ appearing in the SDE (\ref{feller}) bear a different
meaning here with respect to the corresponding parameters $\theta $ and $\mu 
$ in the OU model. Furthermore $\tau <\theta $.

The diffusion coefficient of the process (\ref{feller}) becomes negative if $%
X_{t}<V_{I}$, hence the diffusion interval is $I=[V_I,\infty)$. The boundary
in $V_{I}$ is regular or exit, depeding upon the values of $\mu_2$ and $%
\sigma_2$, according with the Feller classification of boundaries (cf. \cite%
{Karlin}). To determine the transition probability density of the process (%
\ref{feller}) a boundary condition should be added. The natural choice, that
respects the model features, is the reflecting condition:

\begin{equation}
\lim_{x\rightarrow V_I}\left\{\mu(x) f\left( x,t\left\vert
x_{0},t_{0}\right. \right) - \frac{\partial}{\partial x} \left[\sigma^2(x)
f\left( x,t\left\vert x_{0},t_{0}\right. \right)\right]\right\}=0.
\label{reflecting}
\end{equation}

The Feller process is generally known in its standardized form

\begin{equation}
dX_{t}=\left( -\frac{X_{t}}{\tau }+\mu_F\right) dt+\sigma _{2}\sqrt{X_{t}}%
dW_{t} ;\; X_0=-V_I,  \label{feller1}
\end{equation}
that can be easily obtained from (\ref{feller}) by performing the space
transformation $Y_t\rightarrow Y_t - V_I$ and by setting $\mu_F=\mu_2-\frac{%
V_I}{\tau }$. This equation is defined over $I=[0,\infty)$.

A further notation for the parameters of the Feller process, largely
employed in the literature, sets:

\begin{equation}
p \equiv -\frac{1}{\tau};\:q \equiv \mu_F;\:r \equiv \frac{\sigma^2_2}{2}.
\label{pqr}
\end{equation}

The transition pdf of the Feller process $X$ depends (cf. \cite{Karlin})
upon the nature of the lower boundary in $x=0$ for the process (\ref{feller1}%
) or in $x=V_{I}$ for the process (\ref{feller}) and on the selected
boundary condition for the solution of its Kolmogorov equation. If we impose
a zero-flux condition (\ref{reflecting}) in the origin, using the notation (%
\ref{pqr}), we obtain the transition pdf (cf. \cite{giorno}):

\begin{eqnarray}
f_{Fe}\left( x,t\left\vert x_{0},t_{0}\right. \right) &\equiv& \frac{%
\partial P\left( X_{t}<x\left\vert X_{t_{0}}=x_{0}\right. \right) }{\partial
x}=\frac{p\left( \frac{x}{x_{0}}e^{-p(t-t_{0})}\right) ^{\frac{q-r}{2r}}}{%
r(e^{p(t-t_{0})}-1)}  \\
&\times& \exp \left\{ -\frac{p\left( x+x_{0}e^{p(t-t_{0})}\right) }{r\left(
e^{p(t-t_{0})}-1\right) }\right\} I_{\frac{q}{r}-1}\left[ \frac{2p\sqrt{%
xx_{0}e^{p(t-t_{0})}}}{r\left( e^{p(t-t_{0})}-1\right) }\right] .  \notag
\label{ddpFel}
\end{eqnarray}%
Here $I_{\eta }(z)$ indicates the modified Bessel function of the first kind
(cf. \cite{abramowitz}) of parameter $\eta $.

The mean trajectory of the process (\ref{feller}) originated in $X_0=x_0$ is

\begin{equation}
E\left( X_{t}|x_0 \right) =x_{0}e^{-t/\tau }+\mu _{2}\tau \left(
1-e^{-t/\tau }\right).  \label{meanFell}
\end{equation}

while its variance is

\begin{equation}
Var\left( X_{t} \left| x_0 \right.\right) = \tau \sigma^2_2\left(1-e^{-\frac{%
t}{\tau} }\right)\left\{\frac{\mu_2 \tau -V_I}{2}\left(1-e^{\frac{t}{\tau}%
}\right)+\left(x_0-V_I\right)e^{-\frac{t}{\tau} }\right\}.  \label{varFell}
\end{equation}

The FPT pdf of the Feller process cannot be obtained in a closed form but it
can be evaluated by employing the methods described in Sect. \ref%
{Section:Math}. Furthermore its Laplace transform is (cf.\cite{raleigh}):

\begin{equation}
g_{\lambda }\left( S\left\vert x_{0}\right. \right) =\frac{\Phi \left( \frac{%
-\lambda }{p},\frac{q}{r};-\frac{px_{0}}{r}\right) }{\Phi \left( -\frac{%
\lambda }{p},\frac{q}{r};-\frac{pS}{r}\right) }.  \label{Laplace_Fel}
\end{equation}%
Here $\Phi $ denotes the Kummer function (cf. \cite{abramowitz}) and the
notation in (\ref{pqr}) has been employed. The mean firing time, when $%
x_{0}<S<V_{E}$, is (cf. \cite{giorno}):

\begin{equation}
E(T)=\frac{\tau }{\mu _{2}\tau -V_{I}}\left( S-x_{0}+\sum_{n=1}^{\infty }%
\frac{\left( S-V_{I}\right) ^{n+1}-\left( x_{0}-V_{I}\right) ^{n+1}}{\left(
n+1\right) \prod_{i=1}^{n}\left( \mu _{2}\tau -V_{I}+i\tau \sigma
_{2}^{2}/2\right) }\right).  \label{ETFeller}
\end{equation}

If $\tau \mu _{2}>>S$ and $\sigma _{2}$ is suitably small, the mean FPT\ can
be approximated with a formula analoguous to (\ref{ETouap1_0}) for the OU
process. Indeed it holds:

\begin{equation}
E(T)\approx -\tau \ln \left(\frac{S-\mu _{2}\tau}{x_{0}-\mu _{2}\tau}\right).
\label{ETfeap1}
\end{equation}

When the crossing is a rare event, i.e. $x_{0}<<S$ or $\sigma_2$ is small, a
result analogous to (\ref{ETouap2}) can be derived (cf. \cite{gioNoRiasymp}):

\begin{equation}
E(T) \approx \frac{S-V_{I}}{\mu _{2}-\frac{S-V_I}{\tau} -\frac{\sigma
_{2}^{2}}{4}} \Gamma \left( \frac{2\left( \mu _{2}\tau -V_{I}\right)}{\sigma
_{2}^{2}\tau} \right) \left[\frac{\sigma _{2}^{2}\tau}{2\left( S-V_{I}\right)%
}\right]^{\frac{2\left( \mu _{2}\tau -V_{I}\right)}{\sigma _{2}^{2}\tau} }
e^{\left[\frac{2\left( S-V_{I}\right)}{\sigma _{2}^{2}\tau}\right]}.
\label{ETfeap2}
\end{equation}%
Here $\Gamma $\ denotes the Gamma function (cf. \cite{abramowitz}).

The second moment of the FPT has been obtained in \cite{giorno}:

\begin{eqnarray}
E(T^2) &=& \frac{2 \tau E(T)}{\mu_2 \tau -V_I}\left[S-V_I+\sum^{\infty}_{k=1}%
\frac{\left(S-V_I\right)^{k+1}}{(k+1)\prod^{k}_{i=1}\left[%
\mu_2\tau-V_I+\frac{i\tau \sigma^2_2}{2}\right]}\right] \\  \notag
&-& \sum^{\infty}_{k=1}\frac{2 \tau^2\left[(S-V_I)^{k+1}-(x_0-V_I)^{k+1}%
\right]}{(\mu_2 \tau -V_I)(k+1)} \\ \notag
&\times& \sum^{k}_{j=1}\left\{j\prod^{k}_{i=1}\left[%
\mu_2\tau-V_I+ \frac{i\tau\sigma^2_2}{2}\right]\right\}^{-1}.  
\label{ET2Feller}
\end{eqnarray}

The expression for the third moment can be found in \cite{Ricciardi}.

A further variant of Stein's model was proposed in \cite{lanskaetal}
with two reversal potentials, an inhibitory (lower) one $V_I$ and an
excitatory (upper) one $V_E$:

\begin{eqnarray}
dX_t&=& -\frac{1}{\theta}X_t dt +\delta^+(V_E-X_t)dN^+_t+\left[%
\delta^-(X_t-V_I)\right.  \notag \\
&-& \left.J\sqrt{(V_E-X_t)(X_t-V_I)}\right]dN^-_t  \notag \\
X_0&=&x_0.  \label{Steinrev}
\end{eqnarray}
Here the two independent Poisson processes $N^+$ and $N^-$ have intensities $%
\lambda^+$ and $\lambda^-$, respectively, $-1<\delta^-<0<\delta^+<1$, $J$ is a r.v. defined
over the interval $(-1-\delta^-,-\delta^-)$ such that $E(J)=0$. For a sequence
of models (\ref{Steinrev}) indexed by $n$ one can assume that $\delta^+_n,
\delta^-_n\rightarrow 0_-$, $\lambda^+_n, \lambda^-_n\rightarrow \infty$ in
such a way that $\delta^+_n \lambda^+_n \rightarrow \mu\geq 0$, $\delta^-_n
\lambda^-_n \rightarrow \nu \leq 0$. Simultaneously $E(J^2_n)\rightarrow 0_+$
in such a way that $\lambda^-_nE(J^2_n) \rightarrow \sigma^2_3>0$. This
allows to obtain the diffusion approximation

\begin{equation}
dX_{t}=\left( -\frac{X_{t}}{\tau_3 }+\mu _{3}\right) dt+\sigma _{3}\sqrt{%
\left( X_{t}-V_{I}\right) \left( V_{E}-X_{t}\right) }dW_{t} ;\; X_0=x_0
\label{SDEdrev}
\end{equation}
where the new constants $\tau_3$ and $\mu_3=\mu V_E - \nu V_I$ have been
introduced.

The diffusion interval of this process is $[V_I,V_E]$ and its transition pdf
is solution of the Fokker-Planck equation with initial delta condition and
with reflecting conditions of the type (\ref{reflecting}) on both boundaries $V_I$ and $V_E$.

The first two moments of the membrane potential (cf. \cite{lanskaetal}) are: 
\begin{eqnarray}
E(X_{t}\left\vert x_{0}\right. ) &=&x_{0}e^{-t/\tau _{3}}+\mu _{3}\tau
_{3}(1-e^{-t/\tau _{3}})  \label{EXLLS} \\
E\left( X_{t}^{2}\left\vert x_{0}\right. \right) &=&\left(
V_{E}-V_{I}\right) ^{2}\left\{ \frac{\beta \left( 2\beta +\sigma
_{3}^{2}\right) }{\alpha \left( 2\alpha +\sigma _{3}^{2}\right) }+e^{-\alpha
t}\frac{\left( \beta -\alpha \right) \left( 2\beta +\sigma _{3}^{2}\right) }{%
\alpha \left( \alpha +\sigma _{3}^{2}\right) }\right.  \label{secondLLS} \\
&&+e^{-\alpha t-\sigma _{3}^{2}t}\frac{\left( \beta -\alpha \right) \left(
2\beta -2\alpha +\sigma _{3}^{2}\right) }{\left( 2\alpha +\sigma
_{3}^{2}\right) \left( \alpha +\sigma _{3}^{2}\right) }-e^{-\alpha t}\left[
1-\frac{x_{0}-V_{I}}{V_{E}-V_{I}}\right]  \notag \\
&&\left( \frac{2\beta +\sigma _{3}^{2}}{\alpha +\sigma _{3}^{2}}+e^{-\alpha
t-\sigma _{3}^{2}t}\frac{\left( 2\alpha -2\beta +\sigma _{3}^{2}\right) }{%
\alpha +\sigma _{3}^{2}}\right) \left. +e^{-2\alpha t-\sigma _{3}^{2}t}\left[
1-\frac{x_{0}-V_{I}}{V_{E}-V_{I}}\right] ^{2}\right\} .  \notag
\end{eqnarray}%
Here $\alpha =1/\tau _{3}$ and $\beta =\left( -\alpha V_{I}+\mu _{3}\right)
/\left( V_{E}-V_{I}\right) $ and a typo in \cite{lanskaetal} is corrected.
Use of (\ref{EXLLS}) and (\ref{secondLLS}) allows the computation of $%
Var(X_{t}\left\vert x_{0}\right. )$.

The mean firing time through a boundary $S<V_{E}$ is:

\begin{eqnarray}
E(T) &=& \frac{S-x_{0}}{\beta (V_{E}-V_{I})}  \\
&+&\sum_{n=0}^{\infty }\frac{\Gamma (2\beta /\sigma _{3}^{2}+1)\Gamma
(2\alpha /\sigma _{3}^{2}+1)}{\Gamma (2\beta /\sigma _{3}^{2}+n+2)\Gamma
(2\alpha /\sigma _{3}^{2})}\frac{\left( S-V_{I}\right) ^{n+2}-\left(
x_{0}-V_{I}\right) ^{n+2}}{\beta \left( n+2\right) \left( V_{E}-V_{I}\right)
^{n+2}}.  \notag
\label{ETLLS}
\end{eqnarray}

If the boundary crossing is a rare event, a result analoguous to (\ref%
{ETouap2}) and (\ref{ETfeap2}) holds (\cite{lanskaetal}):

\begin{equation}
E(T)\approx \frac{ S-x_{0}}{ \beta \left( V_{E}-V_{I}\right)} \left(1+ \frac{%
\alpha (S+x_0-2V_I)}{(2 \beta+\sigma^2_{3})(V_E-V_I)}. \right)
\label{ETLLSap2}
\end{equation}

When $\mu _{3}\tau_3 > S$, one gets a result analogous to (\ref{ETouap1_0})
and (\ref{ETfeap1}):

\begin{equation}
E(T)\approx -\tau_3 \ln \left(\frac{S-\tau_3 \mu _{3}}{x_{0}-\tau_3 \mu _{3}}%
\right).  \label{ETLLSap1}
\end{equation}

\subsection{Comparison between different LIF models}

\label{Subsection:Comp}

The mathematical complexity of the FPT problem increases with the attempts
to make the models more realistic. However it is desiderable to avoid the
use of complex models when they do not add any improvement with respect to
the simpler ones.

\begin{figure}[htp]
\centering
\includegraphics[width=11cm]{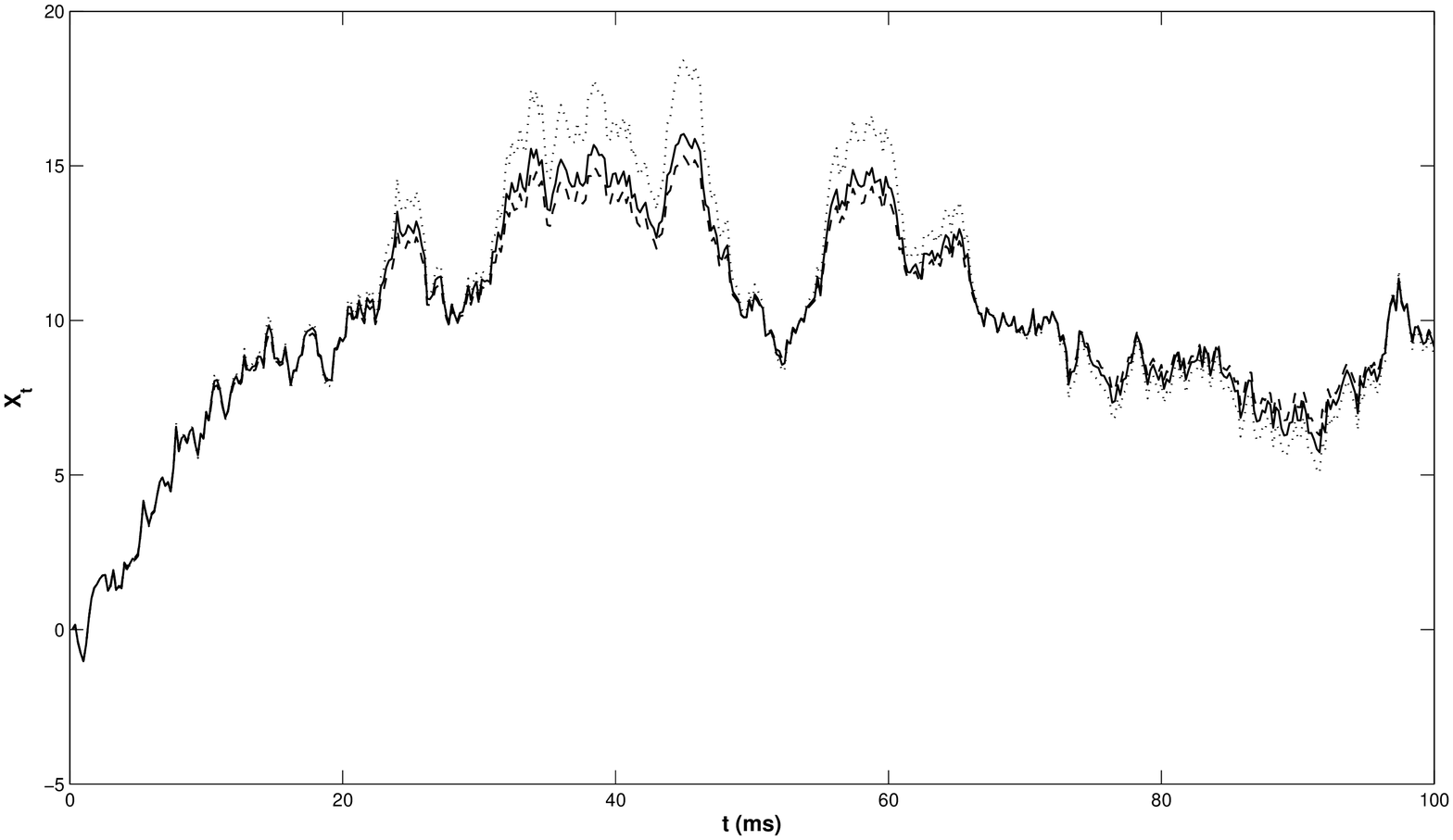}
\caption{Sample paths of the OU (\textit{dashed line}), the Feller (\textit{%
dotted line}) and the double reversal potential (continuous line) models
employing the same leading Wiener process realization. Here $\protect\mu=\protect\mu_2=\protect\mu_3=1$ $mVms^{-1}$, $\sigma^2=0.9$ $mV^2ms^{-1}$, $V_I=-10$ $mV$, $V_E=30$ $mV$.}
\label{fig:Processes}
\end{figure}

In Fig. \ref{fig:Processes} we compare sample paths from the OU, the Feller
process and the process with double reversal potential, simulated using
Euler algorithm (cf. \textbf{SUSANNE}) on the same leading Wiener process trajectory. Furthermore
we consider $\theta =\tau =10$ $ms$ and we choose the same level of
variability at the resting level for all models. Hence $\sigma
_{2}^{2}\left\vert V_{I}\right\vert =\sigma ^{2}$ and $\sigma
_{3}^{2}V_{E}\left\vert V_{I}\right\vert =\sigma ^{2}$. Since the three processes
do not show relevant discrepancies, this analysis suggests to prefer the
first two models due to their better computational tractability.

\begin{figure}[htp]
\centering
\includegraphics[width=14cm]{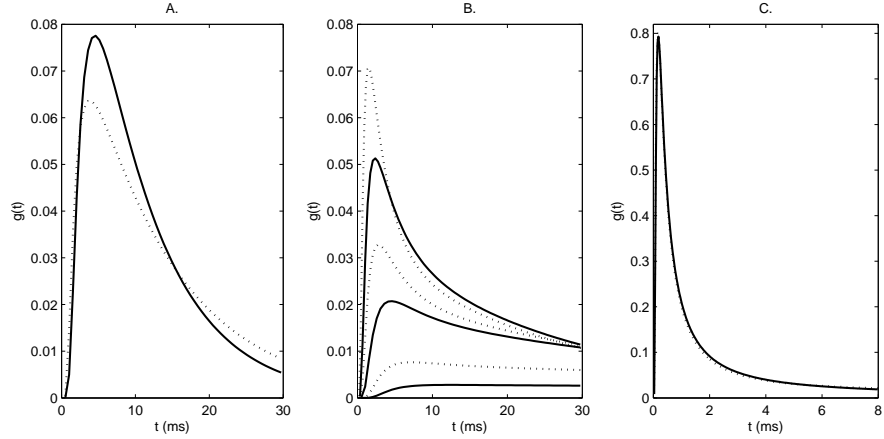}
\caption{Comparisons between the OU (\textit{continuous line}) and the
Feller (\textit{dashed lines}) models. \textbf{Panel A}%
: $V_I=-10$ $mV$, $S=10$ $mV$, $\protect\theta=6.2$ $ms$; parameters
controlling the PSP sizes and the intensities of the input processes: $%
a_e=-i_I=2$ $mV$, $\protect\epsilon=-0.2$, $\protect\lambda=8/\protect\theta%
\protect\cong1.290$ $ms^{-1}$, $\protect\omega=4/\protect\theta\protect\cong%
0.641$ $ms^{-1}$. \textbf{Panel B}: $x_0=0$ $mV$, $V_I=-10$ $mV$, $S=10$ $mV$%
, $\protect\theta=\protect\tau=6.2$ $mV$, $\protect\sigma^2=4,9$ and $16$ 
$mV^2ms^{-1}$ (from bottom to top), $\protect\mu=\protect%
\mu_F=0$ $mVms^{-1}$. \textbf{Panel C}: Feller model $y_0=0$ $mV$, $S=5$ $mV$, $%
V_I=-10$ $mV$, $\protect\tau=6.2$ $ms$, $\protect\mu_2=0$ $%
mVms^{-1}$, $\protect\sigma_2^2=4$ $mVms^{-1}$; OU model $x_0=0$ $mV$, $S=5$ $mV$, 
$\protect\theta=6.2$ $ms$, $\protect\mu=-0.799$ $mVms^{-1}$, $\protect\sigma%
^2=48.03$ $mV^2ms^{-1}$.}
\label{fig:OUFel}
\end{figure}

When one wish to compare the FPT pdf's one gets different results according
with the selected criterium for the parameters values. In \cite{tomasSac2}
the OU and of the Feller ISIs, computed through (\ref{numerical}), are
compared, according with three different criteria:

\begin{itemize}
\item to get the same values for their corresponding discrete versions:

\begin{equation*}
\mu=\mu_2; \;\; \sigma^2 = -\sigma_2^2 V_I,
\end{equation*}

$\theta$ and $\tau$ chosen accordingly (cf. Fig. \ref{fig:OUFel}A);

\item to get the same mean trajectory for both models (Fig. \ref{fig:OUFel}%
B):

\begin{eqnarray*}
\theta &=& \tau;\;\mu = \mu_2 \\
\sigma &=& \sigma_2 \sqrt{-V_I};
\end{eqnarray*}

\item to get almost equal FPT densities. For this aim, one fixes the
parameters for one model and determines a set of parameters,
reproducing a similar ISI distribution, for the other one (cf. Fig. \ref%
{fig:OUFel} C). To guess possible set of parameters for the second process we
impose the equality of the mean and the variance of their FPT's.
\end{itemize}

The last case illustrates the evenience where a histogram of experimentally
obtained ISI's can be fitted by either of the two model distributions.

The use of the more complex Feller model seems preferable when one has data of
the membrane potential between consecutive spikes (\cite{tomasSac2}). When
only the ISI distribution is available, both models fit the data. In (\cite%
{sacsmi}) qualitative comparisons between the OU and the
Feller processes, obtained through the stochastic ordering techniques (cf. 
\cite{sacsmithMCAP}) in Subsubsect. \ref{Subsubsection:Stor}, are presented. In (%
\cite{smithsac}) the same techniques are used for a sensitivity analysis on
the parameters of the FPT pdf. Stochastic ordering properties of the FPT's
are used in \cite{dicrescenzo} to select the model.

Membrane potential data analyzed in \cite{hopfner} show that the same
neuron, under different experimental conditions, can be described either by
the OU model, by the Feller model or by an alternative model with a
quadratic diffusion coefficient.

\subsection{Jump diffusion models}

\label{Subsection:JD}

In Subsections (\ref{Subsection:OU}) and (\ref{Subsection:RevPot}), to
perform the diffusion limits, it was assumed that all the contributions to
the changes in the membrane potential were of the same small amplitude and
the frequency was large enough (cf. \cite{capocelli}). However PSP's
impinging on the soma can play a different role with respect to the
contributions on different areas of the neuron.

LIF models where the subthreshold membrane potential dynamics is described
by jump diffusion processes allow to separate inputs according to their
strongness. Jump diffusion models can be obtained from a variant of the
Stein-type model:

\begin{eqnarray}
dX_{t} &=&-\frac{X_{t}}{\theta }dt+\sum_{j=1}^{n}\delta
_{j}^{+}dN_{t}^{+,j}+\sum_{k=1}^{m}\delta _{k}^{-}dN_{t}^{-,k}+\delta
^{e}dN_{t}^{e}+\delta ^{i}dN_{t}^{i}  \notag \\
X_{t_{0}} &=&x_{0}.  \label{Steinjump}
\end{eqnarray}%
Here $N_{t}^{e},N_{t}^{i}$ are independent Poisson processes of parameters $%
\lambda ^{e}$ and $\lambda ^{i}$ and amplitude $\delta ^{e}$ and $\delta
^{i} $ accounting for the strong contributions. $N_{t}^{+,j},\;N_{t}^{-,k}$
are independent Poisson processes of parameters $\lambda _{j}^{+}$ and $%
\lambda _{k}^{-}$, independent from $N_{t}^{e}$ and $N_{t}^{i}$. If $\delta
_{j}^{+},\delta _{k}^{-}\rightarrow 0$ and at the same time $\lambda
_{j}^{+},\lambda _{k}^{-}\rightarrow \infty $ so that $\delta
_{j}^{+}\lambda _{j}^{+}+\delta _{j}^{-}\lambda _{j}^{-}\rightarrow \mu $
and $(\delta _{j}^{+})^{2}\lambda _{j}^{+}+(\delta _{j}^{-})^{2}\lambda
_{j}^{-}\rightarrow \sigma ^{2}$, a diffusion approximation can be performed
to get a process solution of the SDE:

\begin{equation}
dX_{t}=\left( -\frac{X_{t}}{\theta }+\mu \right) dt+\sigma ^{2}dW_{t}+\delta
^{e}dN_{t}^{e}+\delta ^{i}dN_{t}^{i};\;X_{t_{0}}=x_{0}  \label{jumpDiff}
\end{equation}%
where $W$ is independent from $N_{t}^{e}$ and $N_{t}^{i}$ . The model (%
\ref{jumpDiff}) is a jump diffusion process with an OU underlying
diffusion. Other jump diffusion models may be obtained introducing the
reversal potentials. All these models are of LIF type, requiring the
superposition of a boundary $S$ to mimic the spike activity. The crossings
can occur either for diffusion or for an upward jump when $X\in (S-\delta
^{e},S)$. Hence the spike time is the time of first entrance (FET) into the
strip $(S,\infty )$. The cases of underlying Wiener with drift and OU
process have been considered in (\cite{musilaLansky}, \cite%
{GirSirSacBiosystemPlymouth}, \cite{sacSir}, \cite{sirovich}). The
exponential distribution for the jump epochs preserves the Markov property
of the process (\ref{jumpDiff}). In \cite{sacSir} and \cite{sirsaceBudap}
the case of IG distributed jump epochs is discussed but
in this instance the resulting process is no more a Markov one.

\begin{figure}[htp]
\centering
\includegraphics[width=12cm]{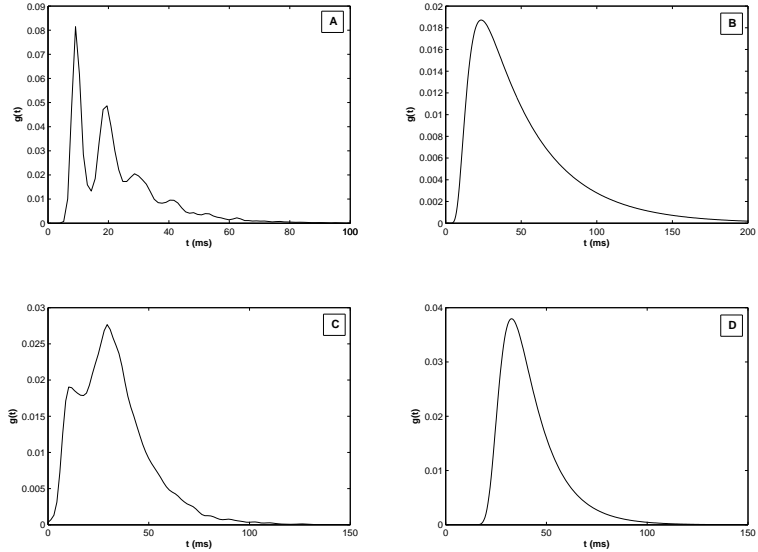}
\caption{Examples of ISI distributions for jump-diffusion processes (\protect
\ref{jumpDiff}); underlying OU diffusion process with different
parameters for Panels A and B, C and D. \textbf{Panel A}: IG time distributed jumps. \textbf{Panel C}: Exponential jump
interarrivals. \textbf{Panel B}, \textbf{Panel D}: ISI distribution for pure
OU diffusion, parameters as in Panels A and C respectively.}
\label{fig:Jumpdiff}
\end{figure}

To compute the ISI distribution for IG and exponential time distributed
jumps one resorts to simulation techniques. Differently from the unimodal
behavior of the ISI distribution of diffusion models, jump diffusion ones
have a multimodal shape (cf. Fig. \ref{fig:Jumpdiff}).

The only analytical results on the FET problem for jump diffusions refer to
an underlying Wiener process with constant boundary. Lower bounds are proposed
in \cite{diCrescetal} for the FET density and in \cite{GirSacSirJ} for the
FET mean and variance, together with exact formulae for the specific case of
large jumps, when the jumps are driven by a counting process. An approximate
solution to an integral equation for the FET density of a jump diffusion
process is discussed in \cite{Wjumps} for the Wiener process.

\subsection{Boundary shapes}

\label{Subsection:Shapes}

\begin{figure}[htp]
\centering
\includegraphics[width=9.cm, height=6cm]{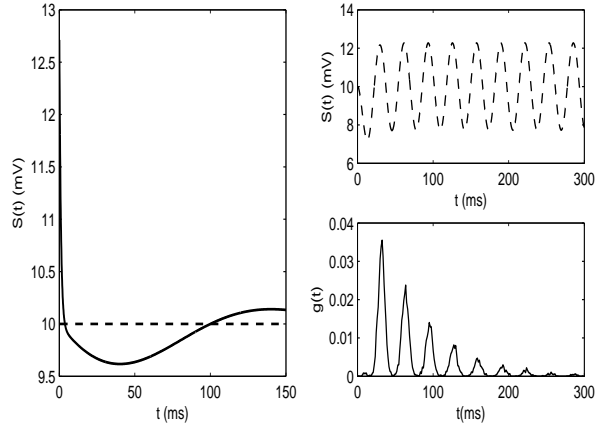}
\caption{\textbf{Panel A:} Time dependent threshold $S(t)=10+3 e^{-t}-0.6
e^{-\frac{t}{100}}\sin(\frac{\protect\pi t}{100})$. \textbf{Panel B:}
Threshold of the transformed process (upper) and simulated FPT pdf
(lower) for an OU process with additional sinusoidal term in the drift
coefficient.}
\label{fig:Tdthresh}
\end{figure}

Constant thresholds are typically employed in LIF models for their easier
mathematical tractability. However the refractory period following each
spike has been modelled by means of threshold shapes (cf. \cite{Ricc}, \cite%
{soglia}). These boundaries attain very high values after a spike, then
decrease under the reference value and finally oscillate around a constant
value (cf. Fig. \ref{fig:Tdthresh}, panel A).

In \cite{longtin3} a dynamic threshold obeying to a differential equation is
considered for the same aim. A boundary which is a linear combination of two
exponentials with different time constants is proposed in \cite{shinomoto},
to fit experimental data. The use of this boundary, together with the lack of
the resetting of the membrane potential after a spike, allows a very goood
fit of the data. A computational method that can reproduce and predict a
variety of spike responses has been deviced in \cite{shinomoto2} using a
multi-timescale adaptive threshold predictor and a nonresetting leaky
integrator.

In \cite{pakLongtin}, \cite{chacron} thresholds with fatigue are proposed to
account for experimental data showing a progressive decrease of excitability
during high frequency firing. This type of threshold destroys the renewal
and Markov character of the process but allows to describe adaptation
phenomena through LIF models. The inclusion of time dependent boundaries
prevents the use of many mathematical methods described in the next Section,
however reliable numerical and simulation techniques can be applied (cf.
Sect. \ref{Section:Math}).

Finally, the study of periodic boundaries or of noisy thresholds (\cite%
{longtin2}) becomes a useful mathematical method to deal with periodic
inputs. Indeed one transforms the original process with time periodic drift
and constant boundary into a time homogeneous diffusion process, constrained
by a periodic absorbing boundary (cf. \cite{transfor1}, \cite{pakSato}, \cite%
{pakdamanMestivier}). To illustrate this idea let us consider an OU model
with periodic input of frequency $\omega $, phase $\varphi $ and amplitude $%
A $. The SDE (\ref{SDE_OU}) for this LIF model is

\begin{equation}
dX_t=\left(-\frac{X_t}{\theta}+\mu+A\sin(\omega t+\varphi)\right)dt+\sigma
dW_t; \; X_0=x_0  \label{OUosc}
\end{equation}
with $x_0<S$. $X$ is not an OU process, however the change of space
variable

\begin{equation}
Y_t=X_t-\frac{A\theta}{\sqrt{1+(\omega\theta)^2}}\left[\sin(\omega
t+\varphi-\xi)-e^{-t/\theta}\sin(\varphi-\xi) \right]  \label{OUmod}
\end{equation}
with $\xi=\arctan(\omega\theta)$ transforms (\ref{OUosc}) into the SDE of an
OU process with parameters $\theta$, $\mu$ and $\sigma$, $%
y_0=x_0$, and the constant boundary $S$ into

\begin{equation}
S(t)=S-\frac{A\theta}{\sqrt{1+(\omega\theta)^2}}\left[\sin(\omega
t+\varphi-\xi)-e^{-t/\theta}\sin(\varphi-\xi) \right].  \label{threshmod}
\end{equation}

The ISI's of the periodically modulated LIF model with constant threshold
are distributed as the ISI's of a LIF model with constant input and
appropriate time-dependent threshold (cf. Fig. \ref{fig:Tdthresh}).

\subsection{Further models}

\label{Subsection:Further}

New efforts on LIF models are mainly devoted to the study of input-output
relationships or to the analysis of small neuronal networks with units
described by LIF models (cf. \cite{sirovich}). New variants of LIF models
have recently appeared in the literature to catch further features such as
plasticity (cf. \cite{gerstner}) or to improve their flexibility and their
predictive power (cf. for example \cite{adex}, \cite{adex2}, \cite{luigia}).
We quote also the LIF compartmental models (cf. \cite{compart1}, \cite%
{compart2}, \cite{compart3}, \cite{compart4}) that extend the one
dimensional ones by introducing systems of SDE's to describe the dynamics of
different components of the neuron such as the dendritic and the soma zone
in the case of two compartmental models. These models require the study of
the FPT of one component through a boundary to describe the ISI's. Up to now
this problem can be handled only through simulations.

\subsection{Refractoriness and return process models}

\label{Subsection:Return}

\begin{figure}[htp]
\centering
\par
\includegraphics[width=10cm, height=6cm]{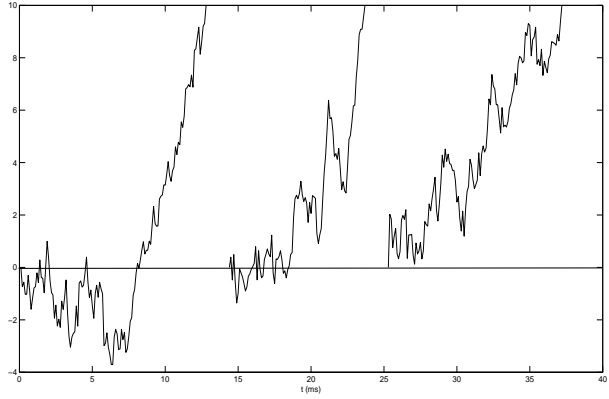}
\caption{Sample path of a return process $Z$ with a constant refractory period.}
\label{fig:Refract}
\end{figure}

An alternative approach to the study of spike trains focus on the number of
spikes in a prescribed time interval. This study allows to introduce the
refractoriness of the neuron in a quite natural way. To this aim one can
associate a return process $\left\{ Z=Z_{t},:t\geq 0\right\}$ in the
interval $(l,S)$, $S\in I$, to any regular diffusion process $X=\left\{
X_{t}:t\geq 0\right\} $ on $I=(l,r)$ as follows. Starting at $x_{0}\in
(l,S) $ at time $0$, the process $Z$ coincides with $X$ until it
attains the level $S$. At this time it is blocked on the boundary $S$ for a
random time and no new crossings can occur during this refractory period.
Then $Z$ and $X$ are instantaneously reset at $x_{0}$ and the
process $Z$ evolves as $X$ until the boundary $S$ is reached again,
and so on. We show in Fig. \ref{fig:Refract} a sample path of this process
when the refractory period is constant. 

Let $F_i$ and $R_i$, $i=0,1,...$, be the r.v.'s describing the time between
the $i$-th reset and the $(i+1)$-th crossing and the $i$-th refractory
period. For time-homogeneous diffusions, the r.v.'s $F_i$ are iid with pdf $%
g(S,t\left|x_0\right.)$. It is also assumed that the r.v.'s $R_i$, 
$i=0,1,...$, are iid with pdf $\varphi(t)$ depending only on the duration of
the refractory period.

A counting process $M=\left\{M_t, \: t \geq 0\right\}$ can be introduced to
describe the number of attainments of the level $S$ by the process $Z$ up to
time $t$. Let

\begin{equation}
q_k(t\left|x_0\right.)=P\left\{M_t=k \left|Z_0=x_0\right. \right\}, \;
k=0,1,...  \label{qkappa}
\end{equation}
be the probability that $k$ firings occur up to $t$. Then (cf. \cite%
{RiccReturn1}):

\begin{eqnarray}  \label{q0kappa}
q_0(t\left|x_0\right.)&=&1-\int^{t}_{0}g(S,\tau\left|x_0\right.)d\tau \\
q_k(t\left|x_0\right.)&=&[g(S,t\left|x_0\right.)*\varphi(t)]^{(k)}*\left[%
1-\int^{t}_{0}g(S,\tau\left|x_0\right.)d\tau\right]+g(S,t\left|x_0\right.) \;
\notag \\
&*& [\varphi(t)*g(S,t\left|x_0\right.)]^{(k-1)}*\left[1-\int^{t}_{0}\varphi(%
\tau)d\tau\right] \;\; k=1,2,...
\end{eqnarray}
where $*$ means convolution and exponent $(k)$ denotes $(k)$-fold
convolution.

Expressions for such probabilities have been obtained in \cite{AlbanoWiener}
for the Wiener process for exponentially distributed refractory periods. In
the general time homogeneous case

\begin{equation}
E\left\{M_{t}^{r}\left\vert x_{0}\right. \right\}
=\sum_{k\geq 1}k^{r}q_{k}(t\left\vert x_{0}\right. )\;\;r=1,2,...
\label{momentiM}
\end{equation}%
is the $r$-th order moment of $M. $Let $I_{i}$, $i=0,1,2,...$ be the r.v.'s
describing the ISI's and let $I_{0}$\ be the time of the first firing. One
has (cf. \cite{RiccReturn2}):

\begin{eqnarray}
E(I)&=&t_1(S\left|x_0\right.)+E(R);\;
E(I^2)=t_2(S\left|x_0\right.)+E(R^2)+2E(R)t_1(S\left|x_0\right.)  \notag \\
E(I^3)&=&t_3(S\left|x_0\right.)+E(R^3)+3E(R)t_2(S\left|x_0%
\right.)+3E(R^2)t_1(S\left|x_0\right.)  \label{momentsISI}
\end{eqnarray}
where $t_r(S\left|x_0\right.)$ is the $r$-th order moment of the FPT. If the
first three moments of the refractory period are finite, the following
asymptotic expressions for large times hold for the first two moments of $M$ (cf. \cite{RiccReturn2}):

\begin{eqnarray}  \label{momretasymp}
E\left\{ M_t \left|x_0\right.\right\} &\cong& \frac{1}{E(I)}t+\frac{1}{2}%
\frac{E(I^2)}{E^2(I)} - \frac{t_1(S\left|x_0\right.)}{E(I)}; \\
E\left\{ \left[M_t\right]^2 \left|x_0\right.\right\} &\cong& \frac{1}{E^2(I)}%
t^2 + \left[ \frac{2E(I^2)}{E^3(I)}-\frac{1}{E(I)}-\frac{2t_1(S\left|x_0%
\right.)}{E^2(I)}\right]t \\
&+& \frac{3}{2}\frac{E^2(I^2)}{E^4(I)} - \frac{2}{3}\frac{E(I^3)}{E^3(I)} - 
\frac{1}{2}\frac{E(I^2)}{E^2(I)} + \frac{t_1(S\left|x_0\right.)}{E(I)}\frac{%
}{}  \notag \\
&+& \frac{t_2(S\left|x_0\right.)}{E^2(I)} - \frac{2E(I^2)}{E^3(I)}%
t_1(S\left|x_0\right.).  \notag
\end{eqnarray}

In \cite{NobileReturn} and \cite{NobileReturn2} the case of absence of
refractory period and that of a random distribution for the return value are
discussed for the Wiener, OU and Feller processes. Alternatively \cite%
{RiccReturn3} proposes to model refractoriness through return processes
characterized by an elastic boundary as firing threshold.

\section{Mathematical methods for First Passage time problem and their
application to the study of neuronal models}

\label{Section:Math}

We update here previous reviews (\cite{Abrahams}, \cite{ricSato} and \cite%
{Ricciardi}) on the methods available up to now to deal with the FPT problem
for stochastic LIF models.

In the case of  diffusion processes, closed form expressions for the
transition pdf's are determined as solutions of the Kolmogorov or the
Fokker-Plank equations of Subsections \ref{Subsection:Wiener}, \ref%
{Subsection:OU} and \ref{Subsection:RevPot}. The Fourier transform is the
typical method to get these solutions.

Closed form solutions for the FPT pdf refers generally to specific time
dependent boundaries. The constant boundary case is known for the Wiener
process (cf. (\ref{FPTWiener})) or for the OU one if $S=0, \mu=0$
and $x_{0}\neq 0.$ The FPT pdf has been determined for the hyperbolic shape
boundary (\ref{hyperbolic}) for the OU process (cf. \cite{Ricc}) or for the
Feller process (cf. \cite{Groups1}) in the case of boundaries corresponding
to symmetries for these processes. Use of the reflection principle allows to
determine the FPT pdf through Daniels' boundary (\cite{Daniels}, \cite%
{Daniels2}, \cite{Lerche}) while a particular FPT pdf is found through a
symmetry based constructive approach in \cite{giornoNobileRicciardi}. The theory of group transformations is used in \cite%
{Groups1} and \cite{Groups2} to determine the transition pdf's of a Feller
process between the origin and a hyperbolic-type boundary and of an OU
process between a lower reflecting and an upper absorbing constant boundary,
but is of no interest for neuronal application and hence we omit the
description of this method.

In Subsect. \ref{Subsection:Analyt} we review the commonly used analytical
techniques for FPT problems: integral equations (\ref{Subsubsection:Int}),
change of variables or of measure (\ref{Subsubsection:Change}), asymptotic
studies (\ref{Subsubsection:Asym}), computation of FPT moments (\ref%
{Subsubsection:Mom}), stochastic ordering (\ref{Subsubsection:Stor}). In
Subsubsect. \ref{Subsubsection:JD} we present the available methods for jump
diffusion processes. In Subsect.(\ref{Subsection:Num}) we introduce the
direct (\ref{Subsubsection:Direct}) and inverse (\ref{Subsubsection:Inv})
FPT problem. Finally in (\ref{Subsection:Sim}) we sketch specific simulation
techniques for FPT's and the numerical tools for their solution.

\subsection{Analytical methods}

\label{Subsection:Analyt}

\subsubsection{Integral equations}

\label{Subsubsection:Int}

In 1943 Fortet (\cite{Fortet}) proved, under mild conditions for the
boundary $S\left( t\right) $, that the Volterra integral equation of the
first kind:

\begin{equation}
f\left(x,t\left\vert x_{0}\right. \right) =\int_{0}^{t}g\left( S( \tau
),\tau \left\vert x_{0}\right. \right)f\left(x,t\left\vert S\left( \tau
\right) ,\tau \right. \right) d\tau,  \label{fortet}
\end{equation}
holding for $x>S(t)$, holds also for $x=S(t)$.

When the boundary is constant and the process is time homogeneous, eq. (\ref%
{fortet}) is of convolution type and the Laplace transform method can be
applied. Denoting as $f_{\lambda}(S\left|x_0\right.)$ and $%
g_{\lambda}(S\left|x_0\right.)$ the Laplace transforms of $%
f(S,t\left|x_0\right.)$ and of $g(S,t\left|x_0\right.)$, for $x>S>x_0$ one
gets:

\begin{equation}
g_{\lambda}(S\left|x_0\right.)=\frac{f_{\lambda}(x\left|x_0\right.)}{%
f_{\lambda}(x\left|S\right.)}.  \label{Laplace_g}
\end{equation}

Generally the Laplace transform cannot be analytically inverted due to its
complex expression (cf. for example (\ref{Laplace})).

Eq. (\ref{fortet}) for $x=S(t)$ has a weakly singular kernel for $\tau
\rightarrow t$. Indeed, any diffusion behaves as a Wiener process for small
times. Hence $f\left( S(t),t\left\vert S\left( \tau \right) ,\tau \right.
\right) \approx \frac{k(t,\tau )}{\sqrt{t-\tau }}$ with $k(t,\tau
)\rightarrow 0$ as $\tau \rightarrow t$. This makes numerical methods for
their solution unstable. Hence it is convenient to switch to a second type
Volterra equation. Integrating (\ref{fortet}) between the left side of the
diffusion interval $l$ and $S(t)$ and then differentiating with respect to
time, one gets a second kind Volterra equation (cf. \cite{RicSacSato}):

\begin{equation}
g(S(t),t\left\vert x_{0}\right. )=2j(S(t),t\left\vert x_{0}\right.
)-2\int_{0}^{t}g\left( S(t),t\left\vert x_{0}\right. \right) j\left( S\left(
t\right) ,t\left\vert S\left( \tau \right) ,\tau \right. \right) d\tau ,
\label{intEqFlux}
\end{equation}%
Here we have introduced the probability current through $z$ at time $u$ of
the diffusion process $X$ whose pdf is solution of (\ref{FP}):

\begin{equation}
j\left( z,u\left\vert w,v\right. \right) =\mu \left( z\right) f\left(
z,u\left\vert w,v\right. \right) -\frac{1}{2}\left. \left\{ \frac{\partial }{%
\partial y}\left[ \sigma ^{2}\left( y\right) f\left( y,u\left\vert
w,v\right. \right) \right] \right\} \right\vert _{y=z}.  \label{flux}
\end{equation}

Eq. \ref{intEqFlux} \ has a weakly singular kernel. For the Wiener process,
when $\left\vert \frac{dS\left( t\right) }{dt}\right\vert \leq Ct^{-\alpha }$%
, with $\alpha <1/2\ $, $C$ a constant and $\lim_{t\rightarrow
0}S(t)>W_{t_0}=x_{0}$, $g\left( S(t),t\left\vert x_{0}\right. \right) $ is
the only $L^{2}$ solution of (\ref{intEqFlux}). It can be expressed as (cf. 
\cite{RicSacSato})

\begin{eqnarray}
g\left( S(t),t\left\vert x_{0}\right. \right) &=& 2j\left( S(t),t\left\vert
x_{0}\right.\right) - 4\int^{t}_{0}d\tau j\left( S(t),t\left\vert S(\tau),
\tau \right.\right) j\left( S(\tau), \tau \left\vert x_{0}\right. \right) 
\notag \\
&+& \sum^{\infty}_{n=1} 4^n \int^{t}_{0}d\tau j_n \left( S(t),t\left\vert
S(\tau), \tau \right.\right) \\
&\times& \left\{2j\left( S(\tau), \tau \left\vert x_{0}\right.\right) - 4
\int^{\tau}_{0}d\theta j\left( S(\tau), \tau\left\vert S(\theta), \theta
\right.\right) j\left( S(\theta), \theta \left\vert
x_{0}\right.\right)\right\}  \notag  \label{infinitesum}
\end{eqnarray}
where

\begin{equation}
j_n\left( S(t),t\left\vert S(\tau), \tau \right.\right) =
\int^{t}_{\tau}d\theta j_1\left( S(\theta),\theta \left\vert S(\tau), \tau
\right.\right) j_{n-1}\left( S(t), t \left\vert S(\theta), \theta
\right.\right)  \label{jenne}
\end{equation}
for $n=2,3,...$ and

\begin{equation}
j_1\left( S(t),t\left\vert S(\tau), \tau \right.\right) =
\int^{t}_{\tau}d\theta j\left( S(\theta),\theta\left\vert S(\tau), \tau
\right.\right) j\left( S(t), t \left\vert S(\theta), \theta \right.\right).
\label{juno}
\end{equation}

A third integral equation was proposed in \cite{buonocore}:

\begin{equation}
g\left( S(t),t\left\vert x_{0}\right. \right) =-2\psi \left(
S(t),t\left\vert x_{0}\right. \right) +2\int_{0}^{t}g\left( S(t),t\left\vert
x_{0}\right. \right) \psi \left( S(t),t\left\vert S\left( \tau \right) ,\tau
\right. \right) d\tau  \label{eqPsi}
\end{equation}
where

\begin{equation}
\psi \left( S(t),t\left\vert x,\tau \right. \right) =\frac{d}{dt}\left\{
F\left( S(t),t\left\vert x,\tau \right. \right) \right\} +k\left( t\right)
f\left( S(t),t\left\vert x,\tau \right. \right) .  \label{psi}
\end{equation}%
Here $F\left( S\left( t\right) ,t\left\vert x,\tau \right. \right)
=\int_{l}^{S\left( t\right) }f\left( y,t\left\vert x,\tau \right. \right) dy$
and $k\left( t\right) $ is an arbitrary continuous function. A suitable
choice for $k\left( t\right) $ allows to make the kernel of (\ref{eqPsi})
regular and hence eq. (\ref{eqPsi}) becomes optimal for numerical
integration. For example the expressions (\ref{psi}) for the OU and the
Feller processes, respectively, that make the kernel of (\ref{eqPsi}%
) regular are (cf. \cite{buonocore})

\begin{eqnarray}
\psi_{OU} \left( S\left( t),t\left\vert x,\tau \right. \right) \right) &=& [ 
\frac{S^{\prime }(t)+S(t)/\theta-\mu}{2}+\frac{e^{\frac{t-\tau}{\theta}}}{%
\theta(1-e^{\frac{2(t-\tau)}{\theta}})} \notag \\   
&\times& \left([S(t)-\mu \theta]e^{\frac{t-\tau}{\theta}}-x+\mu \theta
\right) ] \:f\left( S\left( t),t\left\vert x,\tau \right. \right) \right) 
\label{pou}
\end{eqnarray}
and

\begin{eqnarray}
\psi_{Fel} \left( S\left( t),t\left\vert x,\tau \right. \right) \right) &=&
\frac{p\left[\frac{S(t)^{-p(t-\tau)}}{x} \right]^{\frac{(q-r)}{2r}}}{%
r[e^{p(t-\tau)}-1]}\exp\left\{-\frac{p[S(t)+xe^{p(t-\tau)}]}{%
r[e^{p(t-\tau)}-1]}\right\}  \notag \\
&\times& \left\{\left[S^{\prime }(t)-\frac{pS(t)e^{p(t-\tau)}}{e^{p(t-\tau)}-1}%
+\frac{1}{2}[pS(t)+q-\frac{r}{2}- S^{\prime }(t)]\right]\right. \notag \\ 
&\times& I_{q/r-1}\left[\frac{2p\sqrt{S(t)xe^{p(t-\tau)}}}{r[e^{p(t-\tau)}-1]}%
\right] \notag \\ 
&+& \left.\frac{p\sqrt{S(t)xe^{p(t-\tau)}}}{e^{p(t-\tau)}-1} I_{q/r} %
\left[\frac{2p\sqrt{S(t)xe^{p(t-\tau)}}}{r[e^{p(t-\tau)}-1]} \right]\right\}. 
\label{pfel}
\end{eqnarray}
Other choices for $k(t)$ make the integral on the r.h.s. of (\ref{eqPsi})
vanish for boundaries with particular symmetry properties, such as (\ref%
{hyperbolic}) for the OU process. Infinite sum expansions, bounds and
approximations for $g\left( S(t),t\left\vert x_{0}\right. \right)$ can be
obtained from (\ref{eqPsi}) (cf. \cite{tomasSac}).

Expressions that regularize the kernel of eq. (\ref{eqPsi}) can be
found also in other specific cases.

\subsubsection{Change of variables or change of measures}

\label{Subsubsection:Change}

The transition pdf and the FPT pdf of an assigned processes can be obtained
through changes of variables and/or changes of measure.

Let $Y=\left\{ Y_{t},\;t\geq 0\right\} $ be a diffusion process with $%
I\subseteq \Re $, characterized by drift $\mu (y,t)$ and diffusion
coefficient $\sigma (y,t)$. One may wish to transform this process into a
Wiener process through suitable space-time transformations, when this
transformation exists. In \cite{transfor1} it is shown that a
transformation, conserving the probability mass, mapping the Kolmogorov
equation of the process $Y$ into the analogous equation for the Wiener
process

\begin{equation}
\frac{\partial f^{\prime }}{\partial \tau ^{\prime }}+\frac{\partial
^{2}f^{\prime }}{\partial y^{\prime 2}}=0,  \label{backW}
\end{equation}%
with initial delta condition is of the form

\begin{equation}
\tau ^{\prime }=\phi (\tau );\;y^{\prime }=\psi (y,\tau );\;f(x,t\left\vert
y,\tau \right. )dx=f^{\prime }(x^{\prime },t^{\prime }\left\vert y^{\prime
},\tau ^{\prime }\right. )dx^{\prime }.  \label{transf2}
\end{equation}%
This transformation exists if and only if the infinitesimal moments verify%
\begin{equation}
\mu (y,\tau )=\frac{\sigma _{y}^{\prime 2}(y,\tau )}{2}+\frac{\sigma (y,\tau
)}{2}\left\{ c_{1}(\tau )+\int_{z}^{y}dx\frac{c_{2}(\tau )\sigma ^{2}(x,\tau
)+\sigma _{\tau }^{\prime 2}(x,\tau )}{\left[ \sigma (x,\tau )\right] ^{3}}%
\right\} .  \label{condWiener}
\end{equation}
Here $z\in I$ is an arbitrary value and $c_{1}(t)$ and $c_{2}(t)$ are
arbitrary functions of time. If (\ref{condWiener}) holds the transformation
is:

\begin{eqnarray}
y^{\prime }&=&\psi(y,\tau)=\sqrt{k_1}exp\left[-\frac{1}{2}%
\int^{\tau}_{\tau_0}du c_2(u)\right]\int^{y}_{z}\frac{dx}{\sigma(x,\tau)} 
\notag \\
&-& \frac{\sqrt{k_1}}{2}\int^{\tau}_{\tau_2}du c_1(u)\exp\left[-\frac{1}{2}%
\int^{u}_{\tau_0}d\theta c_2(\theta)\right]+k_2  \notag \\
\tau^{\prime }&=&\phi(\tau)=k_1 \int^{\tau}_{\tau_1}du \exp\left[%
-\int^{u}_{\tau_0}d\theta c_2(\theta)\right]+k_3  \notag \\
f(x,t\left|y,\tau\right.)dx&=&f^{\prime }(x^{\prime },t^{\prime
}\left|y^{\prime },\tau^{\prime }\right.)dx^{\prime }  \label{trWiener}
\end{eqnarray}
where $z\in I$, $\tau_i \in [0,\infty)$ and $k_i$, $i=1,2,3$ are constants
with $k_1>0$. For example the transformation

\begin{eqnarray}
y^{\prime }=\psi(y,\tau)&=& \frac{\sqrt{k_1}}{\sigma}e^{\frac{\tau-\tau_0}{%
\theta}}(y-z)+\frac{\theta \sqrt{k_1}(\frac{z}{\theta}-\mu)}{\sigma}\left[e^{%
\frac{\tau-\tau_0}{\theta}}-e^{\frac{\tau_2-\tau_0}{\theta}}\right]  \notag
\\
&+&k_2  \notag \\
\tau^{\prime }=\phi(\tau)&=& \frac{k_1 \theta}{2} \left[e^{\frac{%
2(\tau-\tau_0)}{\theta}}-e^{\frac{2(\tau_1-\tau_0)}{\theta}}\right]+k_3 
\notag \\
f(x,t\left|y,\tau\right.)dx&=&f^{\prime }(x^{\prime },t^{\prime
}\left|y^{\prime },\tau^{\prime }\right.)dx^{\prime }  \label{OUWiener}
\end{eqnarray}
changes the OU process into a Wiener process. Here $\tau_1,\tau_2>0$ are
arbitrary times. The transformation (\ref{OUWiener}) sends the linear
boundary $S(t)=a+bt$ for the Wiener process into the U-shaped boundary (\ref%
{hyperbolic}) for the OU process.

The Feller process can be transformed into the Wiener one, conserving the
probability mass, only if $\frac{q}{r}=\frac{1}{2}$. In \cite{transfor3} a
necessary and sufficient condition analogous to (\ref{condWiener}) is given
to transform a diffusion process $Y$ into a Feller process.

A change between the measures of two processes is considered in \cite%
{transfor2}, \cite{Alili} and \cite{Alili2} . In \cite{Alili} the Girsanov
theorem (cf. \cite{Roger}) and the change of measure 
\begin{equation}
dP^{OU}=\exp \left[-\frac{1}{2\theta} \left( W^2_t-x^2_0-t \right) - \frac{1%
}{2\theta^2}\int^{t}_{0} W^2_sds \right] dP  \label{Alili}
\end{equation}
are applied to the OU process to obtain its FPT pdf through a constant
boundary. Here $P^{OU}$ and $P$ denote the distributions of an OU process
with $\mu=0$ and of a standard Wiener process, respectively.

\subsubsection{Asymptotic results}

\label{Subsubsection:Asym}

Asymptotic results play an important role in the study of the FPT pdf
because they hold from relatively small values of the involved variable.
Studies on the asymptotic behavior of the FPT pdf belong to two different
classes: large values of the boundary and large
times. Let us first list asymptotic results in the first class. In (\cite%
{nobricSac2}) the asymptotic exponentiality of the FPT for an OU process is
proved; this result is extended in \cite{nobricSac1} to a class of
diffusion processes admitting steady-state density

\begin{equation}
W(x)=\lim_{t\rightarrow\infty}f(x,t\left|x_0\right.) = \frac{c}{\sigma^2(x)}%
\exp\left(\int^x \frac{2\mu(y)}{\sigma^2(y)}dy\right),  \label{steadystate}
\end{equation}
where $c$ is a normalization constant. When the boundary $S$ approaches the
unattainable level $r$ of the diffusion interval, $\lim_{x\rightarrow
r}\sigma^2(x)\left[W(x)\right]^2 E(T)=0$, the following asymptotic result
for the FPT pdf $g(t)$ holds:

\begin{equation}
g(t)\approx \frac{1}{E(T)}\exp\left\{-\frac{t}{E(T)}\right\}.
\label{asympt1}
\end{equation}

Numerical studies on the OU and on the Feller processes show that this
behavior is attained with a negligible error for quite small values of the
boundary $S$ (i.e. for $S=3$ if $\mu=0$, $\theta=1$ and $\sigma^2=1$ for the
OU process). In this asymptotic case the mean FPT, $E(T)$, looses the
dependency upon the initial value $x_0$. Asymptotic results hold for the
same processes in the case of boundaries either asymptotically constant or
asymptotically periodic \cite{gioNoRiasymp}. Periodic boundaries for the OU
process are considered also in \cite{sacer} (see \cite{Ricciardi} for a
review on time dependent boundaries).

Let us now switch to the asymptotic behavior with respect to time. For small
times, the FPT can be approximated with the IG distribution. Indeed near the
origin any diffusion process can be approximated by a Wiener process. In 
\cite{gioNoRiasymp} the asymptotic behaviour, for large $t,$ of the FPT pdf's through
some time-varying boundaries is considered. This paper deals with a class of one
dimensional diffusion processes with steady-state density. The considered
boundaries include periodic boundaries. Sufficient conditions for an
asymptotic exponential behavior are given for the cases of asymptotically
constant and asymptotically periodic boundaries. Explicit expressions are
worked out for the processes that can be obtained from the OU
process by spatial transformations.

The FPT pdf as $t\rightarrow \infty $, for periodic or asymptotically
periodic boundaries $S(t)$, under mild conditions (cf. \cite{gioNoRiasymp})
exhibits damped oscillations with the same period $T$ as the boundary: 
\begin{equation}
g(S(t),t|x_{0})\approx \alpha (t)\exp \left\{ -\int_{0}^{t}\alpha (\tau
)d\tau \right\} .  \label{gasosc}
\end{equation}%
Here $\alpha (t)$ is a periodic function of period $T$:

\begin{eqnarray}  \label{alphaosc}
\alpha(t)&=&-2 \lim_{n\rightarrow \infty}\psi\left[S(t+nT),t+nT
\left|x_0\right. \right] \\ \notag
&=& -\left\{ V^{\prime }(t)+\mu[V(t)]-\frac{d}{dt}\frac{\sigma^2[V(t)]}{4}
\right\}W\left[V(t)\right]
\end{eqnarray}
and $V(t)=\lim_{n\rightarrow \infty}S(t+nT)$. This behavior arises also for
times not too far from the origin. An exponential asymptotic behavior is
also proved for large times and constant boundaries in \cite{Sato_asympt}.

An asymptotic evaluation of the probability that the Wiener process first
crosses a square root boundary is provided in \cite{Sato_asympt_W}. Denoting
as $T_c$ the FPT of the Wiener process trough the boundary $c\sqrt{1+t}$ one
has

\begin{equation}
P(T_{c}>t)\sim _{t\rightarrow \infty }qt^{-p(c)},\;0<p(c)<\frac{1}{2}.
\label{Satoasympt}
\end{equation}%
Here $\lim_{c\rightarrow \infty }p(c)=0;\;\lim_{c\rightarrow 0}p(c)=\frac{1}{%
2}$ and $2p(c)$ is a real solution between $0$ and $1$ with respect to $%
\lambda $ of the equation:

\begin{equation}
\sum^{\infty}_{n=1}\frac{\sin(\frac{\pi \lambda}{2})\Gamma(1+\frac{\lambda}{2%
})(\sqrt{2}c)^n}{\pi n!}\Gamma(\frac{n-\lambda}{2})=1.  \label{Satoasympt3}
\end{equation}

Using the inverse of transformation (\ref{OUWiener}), this result can be
applied to get the asymptotic OU process FPT pdf trough a constant
boundary for large times.

In \cite{SacTom} truncations of the series expansion of the FPT pdf solution
of (\ref{eqPsi}) are used to achieve approximate evaluations. Use of fixed
point theorems is made to obtain bounds for the FPT pdf of the OU and the
Feller processes. Inequalities are proved to find for which times the FPT
pdf can be approximated, within a preassigned error, by means of an assigned
distribution such as the FPT of the Wiener process or the exponential one.

In \cite{nobilePirozzi} the asymptotic behavior of the FPT pdf through
time-varying boundaries is determined for a class of Gauss-Markov processes.

\subsubsection{Moments of the FPT}

\label{Subsubsection:Mom}

Analytical formulae for the moments of the FPT are available only for time
homogeneous processes with constant boundary. Three approaches are possible:

1. derivatives of the Laplace transform of the FPT pdf;

2. solution of second order differential equations;

3. solution of the recursive formula proposed by Siegert (\cite{siegert}).

Having the Laplace transform of the FPT pdf, one can compute:

\begin{equation}
ET^{n}=\left( -1\right) ^{n}\frac{d^{n}g_{\lambda }(S|x_{0})}{d\left(
\lambda \right) ^{n}}  \label{ETn}
\end{equation}%
where $g_{\lambda }(S|x_{0})$ is given by (\ref{Laplace_g}). The presence of
special functions in the Laplace transforms for the OU (\ref{Laplace}) or
for the Feller processes (\ref{Laplace_Fel}) leads to very complex
computations.

Alternatively, using the Kolmogorov equation and eq. (\ref{Laplace_g}), one
can show that the moments of the FPT verify the recursive system
of ordinary differential equations:

\begin{equation}
\sigma^2(x_0)\frac{d^2ET^n(x_0)}{dx_0^2}+\mu(x_0)\frac{dET^n(x_0)}{dx_0}%
=-nET^{n-1}(x_0), \; x_0\in (l,S)
\end{equation}
with boundary conditions:

\begin{equation}
ET^0(l)=0, \; ET^0(S)=1.
\end{equation}

When the process admits steady state distribution, one can write the
solution of (\ref{ETn}) through the Siegert formula (cf. \cite{siegert}):

\begin{equation}
ET^n=t_n(S\left|x_0\right.)=n\int^{S}_{x_0}\frac{2dz}{\sigma^2(z)W(z)}%
\int^{z}_{l}W(y)t_{n-1}(S\left|y\right.)dy.  \label{Siegert}
\end{equation}

Due to the numerical difficulties of these formulae, in \cite%
{LanskySac} approximations are proposed for specific processes (cf.
Subsections \ref{Subsection:OU} and \ref{Subsection:RevPot}) together with
suggestions on the use of each one for specific ranges of the parameters.

\subsubsection{Stochastic ordering}

\label{Subsubsection:Stor}

A further technique for the study of the FPT's is the stochastic comparison
of the FPT's from different models (cf. \cite{sacsmithMCAP}, \cite{smithsac}%
, \cite{sacsmi}). Let us consider the FPT's $T_{1}$ and $T_{2}$ of two
diffusion processes $X_{1}$ and $X_{2}$ over $I_{1}=(l_{1},r_{1})$ and $%
I_{2}=(l_{2},r_{2})$ with drifts $\mu _{i}(x)$, $x=1,2$ and diffusion
coefficients $\sigma _{i}(x)$, $i=1,2$ respectively. Let the two processes $%
Y_{1}$ and $Y_{2}$ be obtained from $X_{1}$ and $X_{2}$ through the
transformation

\begin{equation*}
y_i=g_i(x)=\int^{x}_{l_i} \frac{dz}{\sigma_i(z)}, \: i=1,2.
\end{equation*}
Moreover, let $Y_1$ and $Y_2$ verify the inequalities:

\begin{equation}
\mu_{Y_1}(y) \geq \mu_{Y_2}(y) \:\: \forall y \in [0, g_2(S)]; \; \frac{%
d\mu_{Y_2}(y)}{dy} \leq 0 \:\: \forall y \in [0, g_2(S)]  \label{driftorder}
\end{equation}
and $\sigma^2_1(x) \geq \sigma^2_2(x)$.

For $x_0 \in (\max(l_1, l_2),S)$, $S \in (\max(l_1,l_2), min(r_1, r_2))$, $%
x_0<S$, it holds:

\begin{equation}
T_{X_1}(S|x_0) \leq_{as} T_{X_2}(S|x_0).  \label{Torder}
\end{equation}
In (\ref{Torder}), ''as'' means ''almost surely''. Note that $Y_1$ and $Y_2$
are characterized by unit diffusion coefficien and drift

\begin{equation}
\left.\mu_{Y_{i}}(y)=\frac{1}{\sigma_{i}(x)}(-\frac{1}{4}\frac{d\sigma^2_{i}%
}{dx}+\mu_{i}(x))\right|_{x=g^{-1}_{i}(y)}.  \label{newdrift}
\end{equation}

\subsubsection{Jump diffusion processes}

\label{Subsubsection:JD}

The following integral equation for the FET pdf $\widehat{g}\left( \left.
S,t\right| y,\tau \right)$ of the process $X$ in (\ref{jumpDiff}), defined
over $I=[l,S]$ and originated in the state $y$ at time $\tau$, holds (cf. 
\cite{GirSacJumps2}):

\begin{eqnarray}
\widehat{g}\left( \left. S,t\right\vert y,\tau \right) &=&e^{-\lambda \left(
t-\tau \right) }g\left( S,t\left\vert y,\tau \right. \right) +\int_{\tau
}^{t}du\int_{l}^{S}dze^{-\lambda \left( u-\tau \right) } \\
&\times &\{\lambda ^{e}f^{a}\left( z-a,u\left\vert y,\tau \right. \right)
\bigskip \smallskip +\lambda ^{i}f^{a}\left( z+a,u\left\vert y,\tau \right.
\right) I_{\left( l,S-a\right) }\left( z\right) \}  \notag \\
&\times &\widehat{g}\left( \left. S,t\right\vert z,u\right) +\lambda
^{e}e^{-\lambda \left( t-\tau \right) }\int_{S-a}^{S}dzf^{a}\left(
z,t\left\vert y,\tau \right. \right) .  \notag  \label{densgt}
\end{eqnarray}%
Here $\lambda =\lambda ^{e}+\lambda ^{i}$, $I_{A}(\cdot )$ is the indicator
function of the set $A$, the jump amplitudes are $\delta ^{e}=-\delta
^{i}=a$. Furthermore $g\left( S,t\left\vert y,\tau \right. \right) $ and $%
f^{a}(x,t|y,s)$ are the
FPT and the transition pdf in the presence of the boundary $S$ of the
underlying diffusion process. The following
approximate solution:

\begin{eqnarray}
\widehat{g}\left( \left. S,t\right\vert y,\tau \right) &\approx &e^{-\lambda
\left( t-\tau \right) }g_{t}^{S}\left( y,\tau \right) +\lambda
^{e}e^{-\lambda \left( t-\tau \right) }\int_{S-a}^{S}dzf^{a}\left(
z,t\left\vert y,\tau \right. \right)  \notag \\
&+&\lambda ^{e}e^{-\lambda \left( t-\tau \right) }\int_{\tau
}^{t}du\int_{-\infty }^{S-a}dzf_{a}\left( z,u\left\vert y,\tau \right.
\right) g_{t}^{S-a}\left( z,u\right)  \notag \\
&+&\lambda ^{i}e^{-\lambda \left( t-\tau \right) }\int_{\tau
}^{t}du\int_{-\infty }^{S}dzf^{a}\left( z,u\left\vert y,\tau \right. \right)
g_{t}^{S+a}\left( z,u\right)  \label{g1}
\end{eqnarray}%
holds (cf. \cite{Wjumps}) for a Wiener process with drift $\mu $ and
diffusion coefficient $\sigma $ when $\lambda ^{e}>\lambda ^{i}$, $\lambda
^{e}\ll 1$. Here $g_{t}^{\xi }(z,u)=g(\xi ,t|z,u)$. This approximation can
be interpreted in terms of sample path behavior for the process $X$. For jumps of low frequency, but relevant amplitude with respect to $S$, most of the sample paths cross the boundary either for diffusion without
jumps or for an upward jump when $X_{t}\in \lbrack S-a,S)$ or for diffusion
after at most a single (upward of downward) jump. The possible occurrence of
a higher number of jumps is disregarded. Hence this approximation explains
the first two peaks of the observed multimodal behavior exhibited by the FET
pdf (cf. Fig. \ref{fig:Jumpdiff}).

Some results on the moments of two simplified jump diffusion neuronal models
are discussed in \cite{GirSacJumps2} and \cite{GirSacJumps}. In the
\textquotedblright large jumps\textquotedblright\ model the amplitude of
exponentially time distributed jumps is large enough to determine a crossing of the threshold at each jump.
Assuming that the crossing is a certain event, a recursive relation holds for
the FET moments of order $n\geq 1$ of this model:

\begin{equation}
E[T^n]=\frac{n}{\lambda}\left\{ E[T^{n-1}]+(-1)^n\left[\frac{%
dg^{(n-1)}_{\lambda+\theta}(S\left|x_0\right.)}{d\theta^{n-1}}\right]%
_{\theta=0} \right\}.  \label{largejumps}
\end{equation}
Here $g_{\lambda+\theta}(S\left|x_0\right.)$ is the Laplace transform, of
parameter $\lambda+\theta$, of the pure diffusion FPT, $\lambda$ the
frequency of jumps and the superscript $(m)$ denotes the $m$-th derivative
with respect to $\theta$.

In the 'reset model' exponentially time distributed jumps force the membrane
potential to return instantaneously to its resting level $V_R\equiv x_0$.
Then the dynamics restarts anew till the crossing of the boundary or a new
resetting. This model includes both upward and downward jumps of frequencies 
$\lambda_1$ and $\lambda_2$, whose time epochs are described by means of two
independent Poisson processes. One has:

\begin{eqnarray}  \label{resetmodel}
E[T]&=&\frac{1-g_{\lambda}(S\left|x_0\right.)}{\lambda
g_{\lambda}(S\left|x_0\right.)}; \\
E[T^2]&=&\frac{2E[T]}{g_{\lambda}(S\left|x_0\right.)} + \frac{2}{\left[%
g_{\lambda}(S\left|x_0\right.)\right]^2} \left[ \frac{dg_{\lambda+\theta}(S%
\left|x_0\right.)}{d\theta} \right]_{\theta=0}  \notag
\end{eqnarray}
with the notations as in (\ref{largejumps}) but $\lambda=\lambda_1+\lambda_2$%
.

\subsection{Numerical methods}

\label{Subsection:Num}

\subsubsection{The direct FPT problem}

\label{Subsubsection:Direct}

The direct FPT problem requires to determine the FPT pdf for a given process
assuming the transition pdf and the boundary shape to be known. A large
literature exists on numerical methods to solve the integral equations for
the FPT pdf even for time not homogeneous diffusion processes (cf. \cite{anderssen}, \cite{buonocore}, \cite{durbin}, \cite{Favella}, \cite{gutierrez}, \cite{ricSacSato2}, \cite{Ricciardi}). The one proposed in \cite{buonocore} seems to be the fastest and
most efficient. It consists in discretizing (\ref{eqPsi}) when the function $%
k\left( t\right) $ is chosen to get a regular kernel for the second kind
Volterra equation (cf. (\ref{pou}) or (\ref{pfel}) for the OU and the
Feller processes, respectively). Setting $t=t_{0}+kh$, $k=1,2,...$, $h>0$,
the discretized solution of eq. (\ref{eqPsi}) is

\begin{eqnarray}
\widetilde{g}\left( S\left( t_{0}+h\right) ,t_{0}+h\left\vert
x_{0},t_{0}\right. \right) &=&-2\psi \left( S\left( t_{0}+h\right)
,t_{0}+h\left\vert x_{0},t_{0}\right. \right) ;  \\ \notag
\widetilde{g}\left( S\left( t_{0}+kh\right) ,t_{0}+kh\left\vert
x_{0},t_{0}\right. \right) &=&-2\psi \left( S\left( t_{0}+kh\right)
,t_{0}+kh\left\vert x_{0},t_{0}\right. \right) \\ \notag
&&+2h\sum_{j=1}^{k-1}\widetilde{g}\left( S\left( t_{0}+jh\right)
,t_{0}+jh\left\vert x_{0},t_{0}\right. \right) \\ \notag
&\times& \psi \left( S\left( t_{0}+kh\right) ,t_{0}+kh\left\vert S\left(
t_{0}+jh\right) ,t_{0}+jh\right. \right) \text{ \ \ }\; \\ 
&& k=2,3,...  \notag
\label{numeric}
\end{eqnarray}

\begin{figure}[htp]
\centering
\includegraphics[width=10cm, height=6cm]{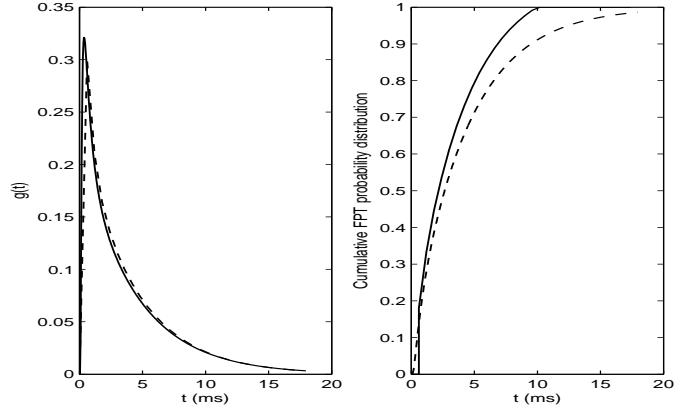}
\caption{FPT pdf (\textbf{panel A}) and FPT probability distribution (%
\textbf{panel B}) of an OU process obtained numerically with $%
\protect\mu=0$ $mVms^{-1}$, $\protect\sigma^2=1$ $mV^2ms^{-1}$, $\protect%
\theta=1$ $ms^{-1}$, $S=1$ $mV$; integration steps $h=0.045$ (%
\textit{continuous line}), $h=0.6$ (\textit{dashed line})}
\label{fig:Intstep}
\end{figure}

A suitable discretization step is necessary to make the numerical
integration reliable. The numerical algorithm (\ref{numeric}) uses previous integration steps to determine the successive values, hence it is
important to get good evaluations on the first intervals. A heuristic rule
is to execute at least twenty integration steps before the maximum of the
FPT pdf occurs. In Fig. \ref{fig:Intstep}, panel A, we show the FPT pdf of a
standard OU process with different integration steps. The error in the FPT
pdf due to a wrong choice for $h$ is enlightened in the evaluation of the FPT
distribution (Panel B).

In \cite{Roman} a strategy is proposed to solve numerically eq. (\ref{eqPsi}%
) with an appropriate choice of the integration step. To this aim a
time-dependent function, the FPT Location function, is introduced.

\subsubsection{The inverse FPT problem}

\label{Subsubsection:Inv}

The inverse FPT problem requires to determine the expression for the
boundary $S(t)$ when the FPT pdf for a diffusion process is known, either in
exact form or through a sample of FPT's. Two numerical algorithms are
proposed to solve this problem for the Wiener process in \cite{zuccaSac} and
extended to the OU process in \cite{saZucca}. The first algorithm proposes a
Monte Carlo procedure to approximate the unknown boundary for the Wiener
process stepwise. This algorithm is reliable and easily implemented but it
is computationally expensive. The second approach, purely numerical, is
computationally more attractive and extensions to processes different from
the Wiener and the OU ones should hold. We present it here in the case of
the Wiener process. By integrating Fortet's equation (\ref{fortet}) in $x$
from $S(t)$ to infinity (cf. \cite{Peskir}) one obtains the integral equation

\begin{equation}
\Psi\!\left(\frac{S(t)}{\sqrt{t}}\right)=\int_0^t\Psi\!\left (\frac{%
S(t)\!-\!S(s)}{\sqrt{t\!-\!s}}\right)g(s) \>ds \hspace{20pt} (t>0)
\label{inteqpesk}
\end{equation}

\noindent where $\Psi (x)=1-\Phi (x)$, $\Phi (x)=\int_{-\infty }^{x}\varphi
(z)\,dz$, $\varphi (z)$ is the standard normal pdf and $g(t)=g(S(t),t\left%
\vert x_{0}\right. )$ is the FPT pdf for the Wiener process through the
threshold $S(t)$. Equation (\ref{inteqpesk}) is a Volterra integral equation
of the first kind in $g(s)$ but it is a non-linear Volterra integral
equation of the second kind in $S(t)$. Its kernel is nonsingular since it is
bounded. Moreover the functions involved in the equation are bounded and $%
\Psi $ is invertible. Hence one can obtain numerically the approximate value 
$S^{\ast }(t_{i})$ of $S(t_{i})$ at the knots $t_{i}=ih$ for $i=1,\ldots ,n$%
; here $h=t/n$ ($t>0$ fixed). The integral on the l.h.s. of (\ref{inteqpesk}%
) is approximated by the Euler method:

\begin{equation}  \label{inteqpesk2}
\Psi\!\left(\frac{S^*(t_i)}{\sqrt{t_i}}\right)=\sum_{j=1}^i\Psi\!\left (%
\frac{S^*(t_i)\!-\!S^*(t_j)}{\sqrt{t_i\!-\!t_j}}\right)g(t_j)\,h \hspace{20pt%
} (i=1,\ldots,n),
\end{equation}

getting a non-linear system of $n$ equations in $n$ unknowns $S^*(t_1),
\ldots,S^*(t_n)$. Note that the $i$-th equation, $i \ge 2$, makes use of the
values $S^*(t_1)$,\ldots,$S^*(t_{i-1})$ in the preceding steps. Hence (\ref%
{inteqpesk2}) can be solved iteratively to get approximate values for the
unknown boundary $S$ at the knots. The convergence of this algorithm and an
estimate of its error are considered in \cite{zuccaSac}.

Extensions to the case in which the FPT pdf is known only through samples of
ISI's make use of the kernel method to approximate the FPT pdf (cf. \cite%
{ZuccaSacerdo}). Applications to neuronal modeling are proposed in \cite%
{saZucca} and \cite{sacVillaZucca} for the OU process. In particular in \cite%
{sacVillaZucca} this algorithm is employed to propose a classification
method of groups of neurons when simultaneously recorded spike trains are
available.

\subsection{Simulation methods}

\label{Subsection:Sim}

In the study of neuronal models, when the process is not time homogeneous or
the boundary is time dependent, the only available general technique is 
simulation. Despite its large use, this methods hides problems that may make
the results for FPT's unreliable. The standard approach to the simulation of
FPT's makes use of discretization algorithms for the SDE describing the
membrane potential dynamics. Various reliable discretization schemes exist
(cf. \cite{Kloeden} and references quoted therein), depending on which
degree of strong or of weak convergence is required. The easiest of these
schemes is the so called Euler-Maruyama one (cf. \textbf{SUSANNE}).

The major cause of error in the simulation of FPT's is related to the chance
to leave the crossing of the boundary undetected due to the discretization
of the sample paths. Indeed a crossing happened inside the discretization
interval results hidden in the discretized sample (cf. Fig. \ref%
{fig:Miscross}). This implies an important overestimation of the FPT, that
does not disappear when the discretization step decreases (cf. \cite%
{girSacSimul}). The number of trials where the error may occur
increases with the decrement of $h$, balancing the corresponding decrease in
the probability of error.

\begin{figure}[htp]
\centering
\includegraphics[width=9cm, height=6cm]{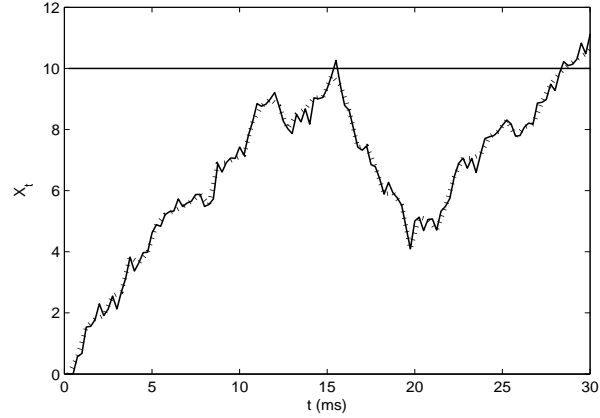}
\caption{A sample path of a diffusion process and its discretization.
Exemplification of a missed crossing of the boundary inside a time interval
of the simulation. \textit{Circles}:
simulated values of the sample path; \textsl{dots}: sample path}
\label{fig:Miscross}
\end{figure}

Different solutions have been proposed to make the simulation of FPT's
reliable. They suggest to evaluate the crossing probability
during each integration step through the bridge process joining the last two
values generated for the diffusion process. A bridge
process $X^{[t_0,t_1]}=\left\{X^{[t_0,t_1]}_t,\:t\in [t_0,t_1]\right\}$ is
associated to a given diffusion $X$ by constraining $X$ to take fixed values
at the time instants $t_0$ and $t_1>t_0$. The process $X^{[t_0,t_1]}$ is
still a diffusion, since its sample paths are a subset of the set of sample
paths of $X$.

The crossing probability of the original diffusion during a simulation step,
of amplitude $h$, coincides with that of its associated bridge on the same
time interval (cf. for example \cite{Roger}). Then one can evaluate, on each
interval, the probability of hidden crossing for this process. For the
Wiener process, the probability that the bridge $W^{[\tau,\tau+h]}$,
originated in the state $y$ at time $\tau$ and constrained to assume the
value $z$ at time $\tau+h$, crosses the boundary $S>y$ during $[\tau,\tau+h]$
is

\begin{equation}
P^*=P_W(S,h,y,z)=\exp\left\{-\frac{2\left(S^2-Sy-Sz+zy\right)}{2h\sigma^2}%
\right\}.  \label{Brbridge}
\end{equation}

One can compare this value with a generated random number $U$ uniform on $%
\left( 0,1\right) $ and, if $U<P^*$, conclude that a crossing has happened
in that interval. The crossing probability of a Wiener process is used to
approximate the crossing of the bridge associated to the process of interest
in \cite{Honerkamp}. To introduce a more precise estimation of the
crossing probability for the bridge we first recall the relationship between
the transition pdf's $f\left( x,t\left| y,\tau \right. \right)$ and $f^{%
\left[ t_{0},t_{1}\right] }\left( x,t\left| y,\tau \right. \right)$ of the
process $X$ and of its bridge $X^{[t_0,t_1]}$ (cf. \cite{Roger}):

\begin{equation}
f^{\left[ t_{0},t_{1}\right] }\left( x,t\left\vert y,\tau \right. \right)
=f\left( x,t\left\vert y,\tau \right. \right) \frac{f\left(
z,t_{1}\left\vert x,t\right. \right) }{f\left( z,t_{1}\left\vert y,\tau
\right. \right) }\text{ \ \ \ \ }t_{0}<\tau <t<t_{1}.  \label{bridgeddp}
\end{equation}

Denoting with $T^{\left[ t_{0},t_{1}\right] }$ the FPT of the bridge $X^{%
\left[ t_{0},t_{1}\right] }$ through the boundary $S$ and with $g^{\left[
t_{0},t_{1}\right] }\left( S,t\left\vert x_{0},t_0\right. \right) $ its pdf
it holds (cf. \cite{Favella}):

\begin{equation}
g^{\left[ t_{0},t_{1}\right] }\left( S,t\left\vert x_{0},t_0\right. \right)
=g\left( S,t\left\vert x_{0},t_0\right. \right) \frac{f\left(
z,t_{1}\left\vert S,t\right. \right) }{f\left( z,t_{1}\left\vert
x_{0},t_{0}\right. \right) }\text{\ \ \ \ }t_{0}<\tau <t<t_{1}.
\label{FPTBridge}
\end{equation}

Integral equations analoguous to (\ref{intEqFlux}) and (\ref{eqPsi}) hold
also for the FPT pdf of the bridge process. In \cite{girSacSimul} an
approximate value of $g^{\left[ t_{0},t_{1}\right] }\left( S,t\left\vert
x_{0}\right. \right) $ is obtained, by disregarding the integral on the
l.h.s of such equations. The approximation using eq. (\ref{eqPsi})

\begin{equation}
g^{\left[ t_{0},t_{1}\right] }\left( S,t\left\vert x_{0},t_0\right. \right)
\cong -2\psi \left( S,t\left\vert x_{0},t_0\right. \right)\frac{f\left(
z,t_{1}\left\vert S,t\right. \right) }{f\left( z,t_{1}\left\vert
x_{0},t_{0}\right. \right) }  \label{apprBridge}
\end{equation}

produces very good estimates in the case of the OU and of the Feller
underlying diffusion processes. In \cite{girSacSimul} the integral of this
approximation over the discretization interval is used to estimate the
probability of a hidden crossing inside each interval. In \cite{gisacZucSim}
a Monte Carlo algorithm is proposed to estimate the crossing probability of
the bridge process. A numerical scheme is applied to the SDE for the bridge
process to generate $N$ samples. If $L$ samples cross the
boundary, the ratio $L/N$ is used to estimate the crossing probability. The
SDE for the bridge has drift and diffusion coefficients:

\begin{equation}
\mu^{\left[ t_{0},t_{1}\right] }(x,t)=\mu(x)+\frac{\sigma^2(x)}{f(z,t_1|x,t)}%
\frac{\partial}{\partial x}f(z,t_1|x,t); \; \sigma^{\left[ t_{0},t_{1}\right]
}(x)=\sigma(x)  \label{SDEbridge}
\end{equation}

respectively (cf. \cite{gisacZucSim}). Here $\mu(x)$ and $\sigma^2(x)$ are
the drift and the diffusion coefficient of the original diffusion $X$.

For a standardized OU process the coefficients (\ref{SDEbridge}) are (cf. 
\cite{gisacZucSim}):

\begin{equation}
\mu^{\left[ t_{0},t_{1}\right] }_{OU}(x,t)= -x+\frac{2\left[ze^{(t_1-t)}-x%
\right]}{\left[e^{2(t_1-t)}-1\right]}; \; \sigma^{\left[ t_{0},t_{1}\right]
}_{OU}(x)= 1,  \label{bridgeOU}
\end{equation}

and for the Feller process are (cf. \cite{gisacZucSim}):

\begin{eqnarray}
\mu^{\left[ t_{0},t_{1}\right] }_{F}(x,t)&=& q + x\left[p-2\frac{%
pe^{p(t_1-t)}}{e^{p(t_1-t)}-1} \right]  \notag \\
&+& \frac{2p\sqrt{xze^{p(t_1-t)}}}{e^{p(t_1-t)}-1} \frac{I_{\frac{q}{r}}%
\left[\frac{2p\sqrt{xze^{p(t_1-t)}}}{r(e^{p(t_1-t)}-1)}\right]}{I_{\frac{q}{r%
}-1}\left[\frac{2p\sqrt{xze^{p(t_1-t)}}}{r(e^{p(t_1-t)}-1)}\right]}  \notag
\\
\sigma^{\left[ t_{0},t_{1}\right] }_{F}(x)&=&\sqrt{2rx}.  \label{bridgeFel}
\end{eqnarray}

Here we employed the notation (\ref{pqr}) and $I(\eta)$ denotes the modified
Bessel function of parameter $\eta$.

A nested algorithm is proposed in \cite{gisacZucSim} for the numerical
solution of the SDE for the bridge. This choice avoids to evaluate the drift
in $t=t_1$ where it becomes singular.

An alternative method to evaluate the hidden crossing probabilities, based
on large deviation arguments, is proposed in \cite{Baldi}. This method is
less precise than the previousluy mentioned ones but it does not request the
knowledge of the transition pdf of the the process $X$.

These algorithms can be applied also to jump diffusion processes but
whenever a jump falls in between the two nodes $t_n$ and $t_{n+1}$ of the
partition, the right end of the time interval $[t_n,t_{n+1}]$ should be
substituted with the epoch $\widetilde{t}_n$, $t_n< \widetilde{t}_n\leq
t_{n+1}$, of the jump event. To account for the possible hidden crossings
inside the discretization intervals, the correction algorithm proposed in 
\cite{girSacSimul} should be employed.

A novel numerical method for the simulation of FPT has been recently
proposed in \cite{magnasco}. The algorithm makes use of the representation
of the stochastic process through an expansion using the Haar functions. It
takes advantage of the dichotomic nature of this development to refine the
description of the process in intervals where posssible hidden crossings may
arise. 

In a recent paper \cite{GirGrSac} it is remarked that the membrane
potential, until the spike, evolves in the presence of the boundary. The SDE for the process constrained by the boundary, i.e. for the
process that has not yet attained the boundary till a fixed time, is determined. The SDE for its bridge, conditioned to cross the boundary for the
first time at its right end, is also determined. Use of the simulation techniques allows to
simulate these stochastic differential equations.

\section{Estimation problems for LIF models}

\label{Section:Est}

A few papers exist on the parameter estimation problem. The literature on
this subject is rather recent, disregarding two oldest papers. The first one
(\cite{LanskyStat}) considers a sample of membrane potential values observed
at discrete times while \cite{inoue} uses the moment method on a sample of
ISI's. We focus here mainly on the available statistical methods for the OU
and the Feller models. Their parameters can be divided into two
groups: the intrinsic parameters, $S,x_{0},V_{I}$ and $\theta$ for the OU
process and $S,x_{0},\tau$ for the Feller process, and the input parameters, 
$\mu$ and $\sigma ^{2}$ for the OU process and $\mu_F$ and $\sigma ^{2}_2$ for
the Feller one. The intrinsic parameters are often disregarded in estimation
problems assuming their direct measure. In \cite{lanskyrefractperiod} the
estimation of the refractory period is also discussed.

We distinguish in the sequel two types of methods, depending upon the available data:

\begin{enumerate}
\item Intracellular membrane recordings;

\item ISI time series.
\end{enumerate}

\subsection{Samples from membrane potential measures}

\label{Subsection:Pot}

In \cite{lanskySanda} a regression method and a maximum likelihood technique
are applied to estimate $\beta =\frac{1}{\theta }$, $\mu $ and $\sigma $
from intracellular data, supposed to follow an OU process. We report here
the maximum likelihood estimators for the case of OU and Feller processes.
One assumes that during an ISI the membrane depolarization $X_{i}=x_{i}$, $%
i=0,1,...,N$, is sampled at the $N+1$ points $t_{i}=ih$. The maximum
likelihood estimates are:

\begin{equation}
\widehat{\beta}=\frac{1}{h}\frac{\sum^{N-1}_{j=0}x^2_j-%
\sum^{N-1}_{j=0}x_{j+1}x_j+(x_N-x_0)\overline{x}}{\sum^{N-1}_{j=0}x^2_j+%
\overline{x}^2N},  \label{betahat}
\end{equation}

\begin{equation}
\widehat{\mu}=\frac{x_N-x_0}{T}+\widehat{\beta}\overline{x},  \label{muhat}
\end{equation}

and 
\begin{equation}
\widehat{\sigma}=\frac{1}{T}\sum^{N-1}_{j=0}\left(x_{j+1}-x_j+x_jh\widehat{%
\beta}-h\widehat{\mu}\right)^2  \label{sigmahat}
\end{equation}

where $\overline{x}=\frac{1}{N}\sum^{N}_{j=0}x_j$, $T=Nh$. The likelihood
estimators $\widehat{\mu}_F$ and $\widehat{\alpha}$ for the parameters $%
\alpha=\frac{1}{\tau}$ and $\mu_F$ of the Feller process (cf. \cite{bibbisor}%
) coincide with (\ref{betahat}) and (\ref{muhat}), while the estimator for $%
\sigma^2_2$, putting $h_{\Delta}=\left(1-e^{-\frac{\Delta}{\tau}}\right)$,
is:

\begin{equation}
\widehat{\sigma}^2_2=\frac{2\sum^{N}_{j=1}X^{-1}_{i-1}\left(X_i-\widehat{\mu}%
_F\tau h_{\Delta}-X_{i-1}e^{-\frac{\Delta}{\tau}}\right)^2}{%
\sum^{N}_{j=1}X^{-1}_{i-1}\tau \left(\widehat{\mu}_F\tau
h^2_{\Delta}+2X_{i-1}h_{\Delta} e^{-\frac{\Delta}{\tau}}\right)}.
\label{sigmahat_Fel}
\end{equation}

Formulae (\ref{betahat}) - (\ref{sigmahat}) are obtained disregarding the
existence of the firing boundary. This approximation determines a bias on
the estimated values (cf. \cite{bibbona}, \cite{bibbona2}). The bias of the
estimator of $\mu $ is of the same order of magnitude as the standard
deviation of the estimate.

Unbiased estimators for $\mu $ are not yet available for the OU and the
Feller models while the bias for the RRW and for the Wiener models are
computed in a closed form in \cite{bibbona}. A comparative study on the
estimators for the Feller process is performed in \cite{bibbona2}. In \cite%
{bibbona3} maximum likelihood estimators are derived from discrete
observations of a Markov process up to the first-hitting time of
a threshold, both in discrete and in continuous time. The models considered
are the RRW, an autoregressive model of order one (AR(1)), and the Wiener,
OU and Feller diffusions. For the last two ones approximations are introduced to evaluate the conditional transition pdf and the FPT
distribution. These approximations hold when the sampling intervals are
small. Their use allow to evaluate the likelihood function.

In \cite{panin} an algorithm is proposed to compute likelihoods, based on
the numerical solution of the integral equation (\ref{eqPsi}). Furthermore
an estimator based on the large deviation principle is suggested to deal
with the case of very small likelihoods in the tails of the distribution.

In \cite{Picchini} a maximum likelihood estimation method is employed for a
particular LIF model with an additional variance parameter modeling possible
slow fluctuations in the parameter $\mu$.

In \cite{hopfner} a sample of discrete observations $%
X_{i\Delta },i_{0}\leq i\leq i_{1},$ $i_{0}:=\left\lceil \frac{t_{0}}{\Delta 
}\right\rceil ,i_{1}:=\left\lceil \frac{t_{1}}{\Delta }\right\rceil $ of the
process $X$ of eq. (\ref{SDE_OU}), over the time interval $[t_{0},t_{1}]$, is considered.
The following nonparametric kernel estimators are proposed: 
\begin{eqnarray}
\widehat{\mu }\left( a\right)  &=&\frac{\sum_{i=i_{0}}^{i_{1}-M}K\left( 
\frac{X_{i\Delta }-a}{h}\right) \left( \frac{X_{\left( i+M\right) \Delta
}-X_{i\Delta }}{\Delta M}\right) }{\sum_{i=i_{0}}^{i_{1}-M}K\left( \frac{%
X_{i\Delta }-a}{h}\right) }  \notag  \label{bhat} \\
\widehat{\sigma ^{2}}\left( a\right)  &=&\frac{\sum_{i=i_{0}}^{i_{1}-M}K%
\left( \frac{X_{i\Delta }-a}{h}\right) \left( \frac{X_{\left( i+M\right)
\Delta }-X_{i\Delta }}{\sqrt{\Delta M}}\right) ^{2}}{%
\sum_{i=i_{0}}^{i_{1}-M}K\left( \frac{X_{i\Delta }-a}{h}\right) }
\label{sigma2hat}
\end{eqnarray}

with $h>0$ a suitable bandwith. Furthermore $K\left( y\right) $ may be
chosen as a rectangular or triangular kernel, for a suitable integer $M$.
Examples of possible kernels are $K(y)=\frac{1}{2}I_{\left\{ -1,+1\right\}
}(y)$ and $K(y)=\left( 1-\left\vert y\right\vert \right) I_{\left\{
-1,+1\right\} }(y)$, with $I_{\left\{ \bullet \right\} }(y)$ indicator
function of the set $\left\{ \bullet \right\} .$

Extensions to other processes and for the selection of the model (OU, Feller
or different ones) are used to exemplify the method. The examples considered
correspond to rarely spiking neurons, a fact that minimizes the problem
underlined in \cite{bibbona}, but prevents its use in other instances.

\subsection{Samples of ISI's}

\label{Subsection:ISI}

The case of ISI data has been recently considered in \cite{MulloIyen}. In
this paper an algorithm is proposed for computing maximum likelihood
estimates with their confidence regions for $\mu $ and $\sigma ^{2}$. The
algorithm numerically inverts the Laplace transform for the OU model. The
method works also to estimate the parameter $\theta $ but it requests larger
samples.

Maximum likelihood estimates for the OU model are obtained in \cite%
{ZhangFeng} using numerical evaluations of (\ref{FPTOUBB}). In \cite%
{ditlevsenLansky1} a variant of the moment method is proposed to estimate
the input parameters of the OU process. In this paper an optimal stopping
theorem is applied to determine the first two exponential moments of the FPT:

\begin{equation}
E\left( e^{T/\theta }\right) =\frac{\mu \theta }{\mu \theta -S},\text{ \ \ \ 
}E\left( e^{2T/\theta }\right) =\frac{2\left( \mu \theta \right) ^{2}-\sigma
^{2}\theta }{2\left( \mu \theta -S\right) ^{2}-\sigma ^{2}\theta }.
\label{expmomOU}
\end{equation}

The moment method is then applied to obtain the
estimators $\widehat{\mu }_{n}$ and $\widehat{\sigma }_{n}^{2}$ from a
sample of ISI's $\left\{ T_{1},...,T_{n}\right\}$:

\begin{equation}
\widehat{\mu }_{n}=\frac{S}{\theta }\frac{Z_{1,n}}{Z_{1,n}-1},\text{ \ \ \ \ 
}\widehat{\sigma }_{n}^{2}=\frac{2S^{2}}{\theta }\frac{Z_{2,n}-Z_{1,n}^{2}}{%
\left( Z_{2,n}-1\right) \left( Z_{1,n}-1\right) ^{2}}  \label{estimOU}
\end{equation}

where 
\begin{equation}
Z_{1,n}=\frac{1}{n}\sum_{i=1}^{n}e^{T_{i}/\theta },\text{ \ \ \ \ }Z_{2,n}=%
\frac{1}{n}\sum_{i=1}^{n}e^{2T_{i}/\theta }.  \label{ZOU}
\end{equation}

This method can be applied only in the suprathreshold region since the
following conditions must be fulfilled: $E\left( e^{T/\theta }\right)
<\infty ,$ $E\left( e^{2T/\theta }\right) <\infty .$ The first condition is
verified if $\mu \theta >S$ and the second holds  if $\mu \theta >S$ and $%
\frac{\sigma ^{2}\theta }{2}<\left( \mu \theta -S\right) ^{2}.$

In \cite{ditlevsenLansky2} analoguous results are proved for the Feller
model, for which

\begin{equation}
E\left( e^{T/\tau }\right) =\frac{\mu_F \tau -y_{0}}{\mu_F \tau -S},\text{ \
\ \ }E\left( e^{2T/\tau }\right) =\frac{\left( \mu_F \tau -y_{0}\right)
^{2}-\sigma ^{2}_2\tau \left( \mu_F \tau /2-y_{0}\right) }{\left( \mu_F \tau
-S\right) ^{2}-\sigma ^{2}_2\tau \left( \mu_F \tau /2-S\right) }.
\label{espMomFeller}
\end{equation}

These expectations converge when $\mu_F \tau >S$ and $\frac{\tau \sigma
^{2}_2}{2}\left( \sqrt{1+\frac{2\mu_F }{\sigma ^{2}_2}}-1\right) <\left(
\mu_F \tau -S\right)$. The estimators are:

\begin{eqnarray}
\widehat{\mu }_{F,n}&=&\frac{S}{\tau }\frac{Z_{1,n}-y_{0}}{Z_{1,n}-1} \\
\widehat{\sigma }_{2,n}^{2}&=&\frac{2\left( S-y_{0}\right) ^{2}}{\tau }\frac{%
\left( Z_{2,n}-Z_{1,n}^{2}\right) }{\left[2\left( Z_{1,n}-1\right) \left(
SZ_{2,n}-y_{0}\right) -\left( SZ_{1,n}-y_{0}\right) \left( Z_{2,n}-1\right)%
\right]}.  \notag  \label{estimFeller}
\end{eqnarray}

Consistency and asymptotic normality of these estimators have been proved in 
\cite{mininni}. The sample sizes required to get the asymptotic conditions
are not huge (some hundreds), hence this
property can be performed in neuronal experiments.

In \cite{ditlevsenditlevsen} an alternative method to estimate $\mu$ and $%
\sigma^2$, based on the analoguous of (\ref{inteqpesk}) for the OU process,

\begin{equation}
\Phi\left(\frac{\mu\theta(1-e^{-\frac{t}{\theta}})-S}{\sqrt{\frac{%
\sigma^2\theta}{2}\left(1-e^{-\frac{2t}{\theta}}\right)}}\right)=%
\int^{t}_{0}g(u)\Phi\left(\sqrt{2}\frac{\mu\theta-S}{\sigma\sqrt{\theta}}%
\sqrt{\frac{1-e^{-\frac{t-u}{\theta}}}{1+e^{-\frac{t-u}{\theta}}}}\right)du,
\label{ditdit}
\end{equation}

is proposed. Numerical results suggest that this approach is preferable to
the previous ones.

The case of observations of the trajectory on very short time intervals is
considered in \cite{lanskyetalshortint}. They propose a method to estimate
the parameters of an OU process in this particular instance.

A recent review (\cite{lanskyDitlevReview}) has appeared summarizing the
state of the art of the estimation problem for diffusion processes in
neuromodeling instances. A comparison of the different estimation methods
for the OU process is performed in \cite{susan}.

Parallel spike trains are deeply discussed in \cite{grun}. This book
presents the methods of correlation analysis together with a review on
different approaches for\ the analysis of single spike trains. The different
chapters discuss many important problems related with the
statistical analysis of spike trains.

Finally the inverse FPT method is applied in \cite{sacVillaZucca} to
classify simultaneously recorded spike trains. The value of the parameters
for the OU process are assumed constant for all the recorded spike trains.
The boundary is determined by the inverse FPT method and a comparison of the
different boundaries is employed to classify the data.

The case of not stationary data is not contemplated by the estimation
procedure. The problem of models whose noise term has a specific temporal
structure has not been solved up to now. In \cite{ZuccaSacerdo} the inverse
FPT method is used on samples of FPT's from an OU process to recover the
boundary shape and to test nonstationary behaviors. The proposed method
makes use of a moving window approach. The inverse FPT\ is applied on
samples from each window. The comparison of the determined boundary
shapes allows to detect changes in the observed dynamics.


\begin{thebibliography}{999}

\bibitem{Abrahams} Abrahams J., A survey of recent progress on
level-crossing problems for random processes. In: Communications and
Networks. A Survey of Recent Advances. (Blake I.F. and Poor H.V., eds.),
6-25. Springer-Verlag, New York (1986)

\bibitem{abramowitz} Abramowitz M. and Stegun I.A., eds. Handbook of
Mathematical Functions with Formulas, Graphs, and Mathematical Tables, New
York, Dover Publications (1972).

\bibitem{AlbanoWiener} Albano G., Giorno V., Nobile A.G. and Ricciardi L.M.
A Wiener-type neuronal model in the presence of exponential refractoriness,
BioSystems 88): 202-215 (2007).

\bibitem{Albano} Albano G., Giorno V., Nobile A.G. and Ricciardi L.M.
Modeling refractoriness for stochastically driven single neurons, Sci. Math.
Jpn. 67(2): 173-189 (2008).

\bibitem{Alili} Alili L., Patie P. and Pedersen J.L. Representation of the
First Hitting Time Density of an Ornstein-Uhlenbeck process, Stochastic
Models 21: 967-980 (2005).

\bibitem{Alili2} Alili L., Patie P. Boundary Crossing Identities fpr
Diffusions having the Time-Inversion Property, J. Theor. Probab. 23: 65--84
(2010).

\bibitem{anderssen} Anderssen R.S., DeHoog F.R. and Weiss R. On the
numerical solution of Brownian motion processes. J. Appl. Prob. 10, 409-418
(1973).

\bibitem{Baldi} Baldi P. and Caramellino L. Asymptotics of hitting
probabilities for general one-dimensional diffusions. Annals of Applied
Probability, 12, 1071-1095 (2002).

\bibitem{bibbona} Bibbona E., L\'{a}nsk\'{y} P., Sacerdote L. and Sirovich
R. Errors in estimation of input signal for integrate and fire neuronal
models. Physical Review E, 78: Art. No. 011918 (2008).

\bibitem{bibbona2} Bibbona E., L\'{a}nsk\'{y} P. and Sirovich R. Estimating
input parameters from intracellular recordings in the Feller neuronal model.
Physical Review E, 81: Art. No. 031916. (2010).

\bibitem{bibbona3} Bibbona E. and Ditlevsen S. Estimation in discretely
observed Markov processes killed at a threshold. Submitted (2010).

\bibitem{bibbisor} Bibby B. and Sorensen M. On estimation for discretely
observed diffusions: A review. Theory Stochastic Process. 2, 49-56 (1996).

\bibitem{adex} Brette R. and Gerstner W. Adaptive Exponential
Integrate-and-Fire Model as an Effective Description of Neuronal Activity,
J. Neurophysiol. 94, 3637-3642 (2005).

\bibitem{buonocore} Buonocore A., Nobile A. G. and Ricciardi L. M. A new
integral equation for the evaluation of the first-passage-time probability
densities. Adv. Appl. Prob. 19, 784-800 (1987).

\bibitem{buonocore2} Buonocore A., Giorno V., Nobile A. G. and Ricciardi L.
M. On the two-boundary first-crossing-time problem for diffusion processes.
J . Appl. Prob. 27, 102-114 (1990).

\bibitem{luigia} Buonocore A., Caputo, L., Pirozzi, E. and Ricciardi L. M.
On a Stochastic Leaky Integrate-and-Fire NeuronalModel. Neural Computation
22, 2558--2585 (2010).

\bibitem{RiccReturn3} Buonocore A., Giorno V., Nobile A.G. and Ricciardi
L.M. A neuronal modeling paradigm in the presence of refractoriness,
BioSystems 67, 35-43 (2002).

\bibitem{Burkitt1} Burkitt A.N., A review of the integrate and fire neuron
model: I. Homogeneous synaptic input. Biol. Cybern. 95, 1-19 (2006).

\bibitem{Burkitt2} Burkitt A.N., A review of the integrate and fire neuron
model: II. Inhomogeneous synaptic input and network properties. Biol.
Cybern. 95, 97-112 (2006).

\bibitem{capocelli} Capocelli R.M. and Ricciardi L.M. Diffusion
approximation and first passage time problem for a model neuron. Kybernetik
8(6):214-23 (1971).

\bibitem{transfor3} Capocelli R.M. and Ricciardi L.M. On the transformation
of diffusion processes into the Feller process, Math. Biosciences 29,
219-234 (1976).

\bibitem{cerbone} Cerbone G., Ricciardi L.M. and Sacerdote L. Mean Variance
and Skewness of first passage time for the Ornstein-Uhlenbeck process. Cyb.
and Systems 12, 395-429 (1981).

\bibitem{longtin3} Chacron M. J., Longtin A., St-Hilaire M. and Maler L.
Suprathreshold stochastic firing dynamics with memory in P-type
electroreceptors. Phys. Rev. Lett. 85, 1576-1579 (2000).

\bibitem{pakLongtin} Chacron M.J,. Pakdaman K,.Longtin A. Interspike
interval correlations, memory, adaptation, and refractoriness in a leaky
integrate-and-fire model with threshold fatigue Neural Computation Volume 15
, Issue 2 253 - 278 (2003).

\bibitem{chacron} Chacron M.J., Lindner B. and Longtin A. Threshold Fatigue
and Information Transmission. Journal of Computational Neuroscience 23:
301-311 (2007).

\bibitem{adex2} Clopath C., Jolivet R., Rauch A., Luscher H.R. and Gerstner
W. Predicting neuronal activity with simple models of the threshold type:
Adaptive exponential integrate-and-fire model with two compartments,
Neurocomput. 70: 1168-1673 (2007).

\bibitem{CoxMiller} Cox D.R. and Miller H.D. The Theory of Stochastic
Processes, Chapman \& Hall (1977).

\bibitem{Daniels} Daniels H. E. The minimum of stationary Markov process
superimposed on a U-shaped trend. J . Appl. Prob. 6, 399-408 (1969).

\bibitem{Daniels2} Daniels H. E. Sequential tests constructed from images.
Ann. Stat., 10 , 394-400 (1982).

\bibitem{dicrescenzo} Di Crescenzo A. and Ricciardi L.M. On a discrimination
problem for a class of stochastic processes with ordered
first-passage-times. Applied Stochastic Models in Business and Industry 17:
205-219 (2001).

\bibitem{diCrescetal} Di Crescenzo A., Di Nardo E. and Ricciardi L.M. On
certain bounds for first-crossing-time probabilities of a jump-diffusion
process. Sci. Math. Jpn. 64, no. 2, 449--460 (2006).

\bibitem{ditlevsenLansky1} Ditlevsen S. and L\'{a}nsk\'{y} P. Estimation of
the input parameters in the Ornstein-Uhlenbeck neuronal model. Physical
Review E 71: Art. No. 011907 (2005).

\bibitem{ditlevsenLansky2} Ditlevsen S. and L\'{a}nsk\'{y} P. Estimation of
the input parameters in the Feller neuronal model Physical Review E 73: Art.
No. 061910 (2006).

\bibitem{ditlevsenditlevsen} Ditlevsen S. and Ditlevsen O. Parameter
estimation from observations of first-passage times of the
Ornstein-Uhlenbeck process and the Feller process. Probabilistic Engineering
Mechanics, 23: 170-179 (2008)

\bibitem{susan} Ditlevsen S. and L\'{a}nsk\'{y} P.: Comparison of
statistical methods for estimation of the input parameters in the
Ornstein-Uhlenbeck neuronal model from first-passage times data. In American
Institute of Physics Proceedings Series, CP1028, Collective Dynamics: Topics
on Competition and Cooperation in the Biosciences, Eds.: L.M. Ricciardi, A.
Buonocore, and E. Pirozzi (2008).

\bibitem{durbin} Durbin J. Boundary crossing probabilities for the Brownian
motion and Poisson processes and techniques for computing the power of the
Kolmogorov Smirnov test J. Appl. Prob. 8, 431-453 (1971).

\bibitem{durbin2} Durbin J. The first-passage density of a continuous
Gaussian process to a general boundary J. Appl. Prob. 22,1, 99-122 (1985).

\bibitem{durbin3} Durbin J. and Williams D. The First-Passage Density of the
Brownian Motion Process to a Curved Boundary J. Appl. Prob. 29,2, 291-304
(1992).

\bibitem{Favella} Favella L., Reineri M.T., Ricciardi L.M. and Sacerdote L.
First-passage-time problems and some related computational methods. Cyb. and
Systems 13, 95-128 (1982)

\bibitem{Fortet} Fortet R. Les fonctions al\'{e}atoires du type Markoff
associ\'{e}es \`{a} certaines \`{e}quations lin\'{e}aires au d\'{e}riv\'{e}%
es partielles du type parabolique J. Math. Pure Appl. (9) 22:177-243 (1943).

\bibitem{gerstein} Gerstein G.L., Mandelbrot B. Random walk models for the
spike activity of a single neuron, Biophys. J. 4,41-68 (1964).

\bibitem{gerstner} Gerstner W. and Kistler W. M. Spiking Neuron Models
Single Neurons, Populations, Plasticity, Cambridge University Press (2002).

\bibitem{raleigh} Giorno V., Nobile A.G., Ricciardi L.M. and Sacerdote, L.
Some remarks on the Raleigh process. J. Appl. Prob. 23, 398-408 (1986).

\bibitem{giorno} Giorno V., L\'{a}nsk\'{y} P., Nobile A.G. and Ricciardi
L.M. Diffusion approximation and first-passage-time problem for a model
neuron: III. A birth-and-death process approach. Biol. Cybern. 58, 387-404
(1988).

\bibitem{giornoNobileRicciardi} Giorno V., Nobile A.G. and Ricciardi L.M. A
symmetry based constructive approach to probability densities for one
dimensional diffusion processes J. Appl. Prob. 27, 707-721 (1989).

\bibitem{giornoNobileRicciardi2} Giorno V., Nobile A.G., Ricciardi L.M. and
Sato S. On the evaluation of first-passage-time probability densities via
non-singular integral equations. Adv. Appl. Prob. 21, 20-36 (1989).

\bibitem{gioNoRiasymp} Giorno V., Nobile A.G. and Ricciardi L.M. On the
asymptotic behavior of first-passage-time densities for one dimensional
diffusion processes and varying boundary. Adv. Appl. Prob. 22, 883-914
(1990).

\bibitem{NobileReturn} Giorno V., Nobile A.G. and Ricciardi L.M.
Instantaneous return process and neuronal firings, in Cybernetics and
Systems Research 1992 (Trappl R. ,Ed.) World Scientific, 829-836 (1992).

\bibitem{RiccReturn2} Giorno V., Nobile A.G. and Ricciardi L.M., C. On the
Moments of Firing Numbers in Diffusion Neuronal Models with Refractoriness.
In: Mira, J., Alvarez, J.R. (Eds.), IWINAC 2005. Lecture Notes in Computer
Sciences 3561, Springer-Verlag, 186-194 (2005).

\bibitem{GirSacJumps2} Giraudo M.T. and Sacerdote L. Some remarks on
First-Passage-Time for Jump-Diffusion Processes. Cybernetics and Systems '96
(R. Trappl Ed.), University of Wien Press, Wien, p. 518-523 (1996).

\bibitem{GirSacJumps} Giraudo M.T. and Sacerdote L. Jump-diffusion processes
as models for neuronal activity, BioSystems 40: 75-82 (1997).

\bibitem{soglia} Giraudo M.T. and Sacerdote L., Simulation methods in
neuronal modeling. BioSystems 48:77--83 (1998).

\bibitem{girSacSimul} Giraudo M.T. and Sacerdote L. An improved technique
for the simulation of first passage times for diffusion processes
Communication in Statistics: simulation and computation. 28, 4 (1999).

\bibitem{Groups2} Giraudo M.T., A similarity solution for the
Ornstein-Uhlenbeck diffusion process constrained by a reflecting and an
absorbing boundary. Ricerche di Matematica. 49(1):47--63 (2000).

\bibitem{gisacZucSim} Giraudo M.T., Sacerdote L. and Zucca C. Evaluation of
first passage times of diffusion processes through boundaries by means of a
totally simulative algorithm. Meth. Comp. Appl. Prob. 3, 215-231 (2001).

\bibitem{GirSirSacBiosystemPlymouth} Giraudo M.T., Sacerdote L. and Sirovich
R. Effects of random jumps on a very simple neuronal diffusion model,
BioSystems 67: 75-83 (2002).

\bibitem{GirSacSirJ} Giraudo M.T., Sacerdote L. and Sirovich, R. Effects of
random jumps on a very simple neuronal diffusion model, BioSystems 67: 75-83
(2002).

\bibitem{mininni} Giraudo M.T., Mininni R. and Sacerdote L. On the
asymptotic behavior of the parameter estimators for some diffusion
processes: application to neuronal models, Ricerche di Matematica 58 (1):
103-127 (2009).

\bibitem{Wjumps} Giraudo M.T. An approximate formula for the
first-crossing-time density of a Wiener process perturbed by random jumps,
Statistics and Probability Letters 79: 1559-1567 (2009).

\bibitem{GirGrSac} Giraudo M.T., Greenwood, P.E. and Sacerdote L. How sample
paths of Leaky Integrate and Fire models are influenced by the presence of a
firing threshold. In press on Neural Computation (2011).

\bibitem{grun} Gr$\ddot{u}$n S. and Rotter P. Analysis of parallel spike
trains. Springer, New York (2010).

\bibitem{gutierrez} Gutierrez, R., Ricciardi, L.M., Rom\'{a}n, P. and
Torres, F., First-passage-time densities for time-non-homogeneous diffusion
processes. J. Appl. Prob. 34, 623-631 (1997).

\bibitem{lanskyrefractperiod} Hampel D, L\'{a}nsk\'{y} P: On the estimation
of refractory period. Journal of Neuroscience Methods, 171: 288-295 (2008).

\bibitem{diesmann} Helias M., Deger M., Diesmann M. and Rotter S.:
Equilibrium and response properties of the integrate-and-fire neuron in
discrete time. Frontiers in Computational Neuroscience, 3, Article 29 (2010).

\bibitem{HH} Hodgkin A. and Huxley A. A quantitative description of membrane
current and its application to conduction and excitation in nerve. J.
Physiol. 117:500-544 (1952).

\bibitem{Honerkamp} Honerkamp J. Stochastic Dynamical Systems: Concepts,
Numerical methods, Data Analysis, VCH (1994).

\bibitem{hopfner} Hopfner R. On a set of data for the membrane potential in
a neuron Math. Biosci., 207(2):275--301 (2007).

\bibitem{inoue} Inoue J, Sato S and Ricciardi LM. On the parameter
estimation for diffusion models of single neuron's activities. I.
Application to spontaneous activities of mesencephalic reticular formation
cells in sleep and waking states. Biol Cybern 1995;73(3):209-221 (1995).

\bibitem{jolivet} Jolivet R., Lewis T.J., and Gerstner W., Generalized
Integrate-and-Fire Models of Neuronal Activity. Approximate Spike Trains of
a Detailed Model to a High Degree of Accuracy, J. Neurophysiol.,
92(2):959--976 (2004).

\bibitem{BiolCybRev} Jolivet R., Roth A., Schurmann F., Gerstner W. and Senn
W. Special issue on quantitative neuron modeling. Biol. Cybern. 99, 237-239
(2006).

\bibitem{Kallianpur} Kallianpur G. On the diffusion approximation to a
discontinuous model for a single neuron. Contributions to statistics,
247-258, North-Holland, Amsterdam, (1983).

\bibitem{Karlin} Karlin S. and Taylor H. M. A Second Course in Stochastic
Processes. Academic Press (1981).

\bibitem{kistler} Kistler W. M., Gerstner W. and vanHemmen J. L. Reduction
of the Hodgkin-Huxley Equations to a Single-Variable Threshold Model, Neural
Comput., 9(5):1015--1045, (1997).

\bibitem{Kloeden} Kloeden P. and Platen P. The numerical solution of
Stochastic differential equations. Springer (1992).

\bibitem{shinomoto} Kobayashi R., Tsubo Y. and Shinomoto S. Predicting spike
times of any cortical neuron. Frontiers in Systems Neuroscience. Conference
Abstract: Computational and systems neuroscience. doi:
10.3389/conf.neuro.06.2009.03.196 (2009).

\bibitem{shinomoto2} Kobayashi R., Tsubo Y. and Shinomoto S. Made-to-order
spiking neuron model equipped with a multi-timescale adaptive threshold.
Frontiers in Computational Neuroscience, 3, Article 9 (2009).

\bibitem{LanskyStat} L\'{a}nsk\'{y} P. Inference for the diffusion models of
neuronal activity. Mathematical Biosciences 67:247-260 (1983).

\bibitem{Lanska} L\'{a}nsk\'{y} P. and L\'{a}nsk\'{a} V. Diffusion
approximations of the neuronal model with synaptic reversal potentials Biol
Cybern. 56 , 19-26 (1987).

\bibitem{musilaLansky} L\'{a}nsk\'{y} P. and Musila M. Generalized Stein's
model for anatomically complex neurons Biosystems 25, 179-191 (1991).

\bibitem{lanskaetal} L\'{a}nsk\'{a} V., L\'{a}nsk\'{y} P. and Smith C.E.
Synaptic transmission in a diffusion model for neural activity. J Theor Biol
166:393-406 (1994).

\bibitem{tomasSac2} L\'{a}nsk\'{y} P., Sacerdote L. and Tomassetti F. (1996)
On the comparison of Feller and Ornstein-Uhlenbeck models for neural
activity. Biol. Cybern. 73, 457-465 (1995).

\bibitem{Lanskysato} L\'{a}nsk\'{y} P. and Sato S. The stochastic diffusion
models of nerve membrane depolarization and interspike interval generation.
Journal of Peripheral Nervous Systems 4: 27-42 (1999).

\bibitem{compart1} L\'{a}nsk\'{y} P. and Rodriguez R. Coding range of a
two-compartmental model of a neuron. Biol. Cybern. 81, 161 (1999).

\bibitem{LanskySac} L\'{a}nsk\'{y} P. and Sacerdote L. The
Ornstein-Uhlenbeck neuronal model with signal-dependent noise. Physics
Letters A 285, 132-140 (2001).

\bibitem{lanskySanda} L\'{a}nsk\'{y} P., Sanda P., He J.F.: The parameters
of the stochastic leaky integrate-and-fire neuronal model. Journal of
Computational Neuroscience 21: 211-223 (2006).

\bibitem{lanskyDitlevReview} L\'{a}nsk\'{y} P. and Ditlevsen S. A review of
the methods for signal estimation in stochastic diffusion leaky
integrate-and-fire neuronal models. Biological Cybernetics 99, 253-262
(2008).

\bibitem{Lapique} Lapique, L. Reserches quantitatives sur l'excitation \'{e}%
lectrique des nerfs trait\'{e}e comme une polarization. J. Physiol. Pathol.
Gen. 9, pp. 620-635 (1907).

\bibitem{Lerche} Lerche, H.R. Boundary Crossing of Brownian Motion, Lecture
Notes in Statistics, Vol. 40, Springer-Verlag, Heidelberg-New York (1986).

\bibitem{longtin2} Lindner B., Chacron M.J. and Longtin A. Integrate and
fire neurons with threshold noise: A tractable model of how interspike
interval correlations affect neuronal signal transmission Physical Review E
72, 021911 (2005).

\bibitem{MulloIyen} Mullowney P. and Iyengar Maximum Likelihood Estimation
of an Integrate and Fire Neuronal Model. Neural Computation.
2008;20:1776-1795 (2007).

Neurocomputing 32-33: 219-224 (2000).

\bibitem{nobricSac1} Nobile A.G., Ricciardi L.M. and Sacerdote L.
Exponential trends of first passage time densities for a class of diffusion
processes with steady-state distribution. J. Appl. Prob. 22, 611-618 (1985).

\bibitem{nobricSac2} Nobile A.G., Ricciardi L.M. and Sacerdote L.
Exponential trends of Ornstein-Uhlenbeck first-passage-time densities. J.
Appl. Prob. 22, 360-369 (1985).

\bibitem{nobilePirozzi} Nobile A.G., Pirozzi E. and Ricciardi L.M. On the
Estimation of First-Passage Time Densities for a Class of Gauss-Markov
Processes EUROCAST 2007,M. Diaz ed. LNCS 4739, pp. 146-153 (2007).

\bibitem{pakdamanMestivier} Pakdaman K. and Mestivier D. (2001) External
noise synchronizes forced oscillators Phys. Rev. E 64, 030901

\bibitem{lanskyetalshortint} Pawlas Z., Klebanov L.B., Prokop M. and L\'{a}%
nsk\'{y} P. Parameters of Spike Trains Observed in a Short Time Window.
Neural Computation, 20: 1325-1343 (2008).

\bibitem{panin} Paninski L., Haith A. and Szirtes G. Integral equation
methods for computing likelihoods and their derivatives in the stochastic
integrate and fire model. J Comput Neurosci (2008) 24:69-79 (2008).

\bibitem{Peskir} Peskir, G. Limit at zero of the Brownian first-passage
density. Probab. Theory Related Fields 124, 100-111 (2002).

\bibitem{Picchini} Picchini U., L\'{a}nsk\'{y} P., De Gaetano A. and
Ditlevsen S. Parameters of the diffusion leaky integrate-and fire neuronal
model for a slowly fluctuating signal, Neural Comput. 20: 2696-2714 (2008).

\bibitem{plesser2} Plesser H.E. and Diesmann M. Simplicity and efficiency of
integrate-and-fire neuron models. Neural Computation, 21, 353-359 (2009).

\bibitem{transfor1} Ricciardi L.M. On the transformation of Diffusion
Processes into the Wiener process. Journal of Mathematical Analysis and
Applications 54, 1: 185-199 (1976).

\bibitem{Ricc} Ricciardi L.M. Diffusion Processes and Related Topics in
Biology, Lecture Notes in Biomathematics, Vol. 14. Springer Verlag, Berlin
(1977).

\bibitem{RicSac} Ricciardi L.M. and Sacerdote L. The Ornstein-Uhlenbeck
Process as a Model for Neuronal Activity, Biol. Cybern. 35, 1-9 (1979).

\bibitem{ricSacSato2} Ricciardi L.M., Sacerdote L. and Sato S. Diffusion
approximation and first passage time problem for a model neuron II Outline
of a computational method. Math. Biosciences, 64, 29-44 (1983).

\bibitem{RicSacSato} Ricciardi L.M., Sacerdote L. and Sato S. On an Integral
Equation for First-Passage-Time Probability Densities Journal of Applied
Probability, Vol. 21, No. 2. pp. 302-314 (1984).

\bibitem{ricSato} Ricciardi L.M. and Sato S. Diffusion processes and
first-passage-time problems Lectures Notes in Biomathematics and Informatics
Ricciardi L.M. ed. Manchester Univ. Press (1989).

\bibitem{NobileReturn2} Ricciardi L.M., Di Crescenzo A., Giorno, V. and
Nobile A.G. On the instantaneous return process for neuronal diffusion
models. In: Marinaro, M., Scarpetta, G. (Eds.), Structure: from Physics to
General Systems. World Scientific, pp. 78-94 (1992).

\bibitem{Ricciardi} Ricciardi L.M., Di Crescenzo A., Giorno V. and Nobile A.
An outline of theoretical and algorithmic approaches to first passage time
problems with applications to biological modeling. Mathematica Japonica 50
(2):247-321 (1999).

\bibitem{RiccReturn1} Ricciardi L.M., Esposito G., Giorno V. and Valerio C.
Modeling neuronal firing in the presence of refractoriness. In: Mira, J.,
Alvarez, J.R. (Eds.), IWANN 2003. Lecture Notes in Computer Sciences 2686,
Springer-Verlag, 1-8 (2003).

\bibitem{compart2} Rodriguez R. and L\'{a}nsk\'{y} P. Two-compartment
stochastic model of a neuron with periodic input. Lecture Notes in Computer
Science 1606 "Foundations and Tools for Neural Modeling (IWANN'99)",
Springer Verlag (1999).

\bibitem{compart3} Rodriguez R. and L\'{a}nsk\'{y} P. A simple stochastic
model of spatially complex neurons. Biosystems 58, 49 (2000).

\bibitem{compart4} Rodriguez R. and L\'{a}nsk\'{y} P. Effect of spatial
extension on noise-enhanced phase locking in a leaky integrate-and-fire
model of a neuron. Physical Review E 62, 8427 (2001).

\bibitem{Roger} Roger L.C.G. and Williams D. Diffusions, Markov Processes
and Martingales. Wiley Series in Probability and Mathematical Statistics
(1987).

\bibitem{Roman} Rom\'{a}n P., Serrano J.J. and Torres F., First-passage-time
location function: Application to determine first-passage-time densities in
diffusion processes. Computational Statistics and Data Analysis 52,
4132--4146 (2008).

\bibitem{sacer} Sacerdote L. Asymptotic behavior of Ornstein-Uhlenbeck
first-passage-time density through boundaries. Applied Stochastic Models and
Data Analysis. 6, 53-57 (1988).

\bibitem{Groups1} Sacerdote L., On the solution of the Fokker-Planck
equation for a Feller process. Adv. Appl. Prob. 22(1):101--110 (1990).

\bibitem{transfor2} Sacerdote L. and Ricciardi L.M. On the transformation of
diffusion equations and boundaries into the Kolmogorov equation for the
Wiener process. Ricerche di Matematica 41, 1: 123-135 (1992).

\bibitem{tomasSac} Sacerdote L. and Tomassetti F. On Evaluations and
Asymptotic Approximations of First-Passage-Time Probabilities. Adva
Applied Probability, Vol. 28, No. 1. pp. 270-284 (1996).

\bibitem{smithsac} Sacerdote L. and Smith C.E. (2000) New Parameter
Relationships Determined Via Stochastic Ordering for Spike Activity in a
Reversal Potential Model. BioSystem 58: 59-65.

\bibitem{sacsmi} Sacerdote L. and Smith C.E. (2000) A qualitative comparison
of some diffusion models for neural activity via stochastic ordering. Biol.
Cybernetics 83, 6 543-551.

\bibitem{sacSir} Sacerdote L. and Sirovich R. Multimodality of the
interspike interval distribution in a simple jump-diffusion model Scientiae
Mathematicae Japonicae Online, Vol. 8, 359-374 (2003).

\bibitem{saZucca} Sacerdote L. and Zucca C. Threshold shape corresponding to
a Gamma firing distribution in an Ornstein-Uhlenbeck neuronal model.
Scientiae Mathematicae Japonicae. 58-2, 295-30 (2003).

\bibitem{ZuccaSacerdo} Sacerdote L. and Zucca C. On the relationship between
interspikes interval distribution and boundary shape in the
Ornstein-Uhlenbeck neuronal model. ECMTB Proceedings. V. Capasso ed.,
161-168 (2003).

\bibitem{sacsmithMCAP} Sacerdote L. and Smith C.E. Almost sure comparisons
for first passage times of diffusion processes through boundaries.
Methodology and Computing in Applied Probability 6, Number 3, 323-341 (2004).

\bibitem{sirsaceBudap} Sacerdote L. and Sirovich R. Noise induced phenomena
in jump-diffusion models for single neuron spike activity. IJCNN
Proceedings, Budapest (2004).

\bibitem{sacVillaZucca} Sacerdote L., Villa A.E.P. and Zucca C. On the
classification of experimental data modeled via a stochastic leaky integrate
and fire model through boundary values. Bull. Math. Biol. 68(6):1257-74
(2006).

\bibitem{Sato_asympt_W} Sato S. Evaluation of the First-Passage Time
Probability to a Square Root Boundary for the Wiener Process, J. Appl. Prob.
14 (4): 850-856 (1977).

\bibitem{Sato_asympt} Sato S. Note on the Ornstein-Uhlenbeck process model
for stochastic activity of a single neuron, Lecture Notes in Biomath, 70:
146-156 (1987).

\bibitem{Segundo} Segundo J., Vibert J.-F., Pakdaman K., Stiber M. and Diez
Martinez O. Noise and the neuroscience: a long history, a recent revival and
some theory. In: Pribram, KH eds. , Origins: brain \& self organization,
Erlbaum, Hillsdale, NJ (1994).

\bibitem{pakSato} Shimokawa T., Pakdaman K. and Sato S. Time-scale matching
in the response of a leaky integrate-and-fire neuron model to periodic
stimulus with additive noise, Phys. Rev. E 59, 3427 3443 (1999).


\bibitem{siegert} Siegert A.J.F. On the First Passage Time Probability
Problem, Phys. Rev. 81, 617--623 (1951).

\bibitem{sirovich} Sirovich R. Mathematical models for the study of
synchronization phenomena in neuronal networks, Ph. D. Thesis, University of
Torino and Universit\'{e} de Grenoble (2006).

\bibitem{stein} Stein R.B. A theoretical analysis of neuronal variability,
Biophys. Journal 5, 385-386 (1965).

\bibitem{magnasco} Taillefumier T. and Magnasco M.O. A Fast Algorithm for
the First-Passage Times of Gauss-Markov Processes with H\"{o}lder Continuous
Boundaries, J. Stat. Phys. 140: 1130--1156 (2010).

\bibitem{Tuckwell1} Tuckwell H. C. Introduction to theoretical neurobiology.
Volume 1. Linear cable theory and dendritic structure, Cambridge University
Press, Cambridge, UK (1988).

\bibitem{Tuckwell2} Tuckwell, H. C. Introduction to theoretical
neurobiology. Volume 2. Introduction to theoretical neurobiology, Cambridge,
UK:Cambridge University Press (1988).

\bibitem{OrnUhl} Uhlenbeck G.E. and Ornstein L.S. On the theory of Brownian
motion, Phys. Rev. 36, 823-841 (1930).

\bibitem{WangPotz} Wang L. and Potzelberger K. Boundary crossing probability
for Brownian motion and general boundaries, J. Appl. Prob. 34: 54-65 (1997).

\bibitem{ZhangFeng} Zhang X., You G., Chen T. and Feng J. K. Maximum
likelihood decoding of neuronal inputs from an interspike interval
distribution, Neural Computation 19 (4), 1319-1346 (2009).

\bibitem{zuccaSac} Zucca C. and Sacerdote L. On the Inverse
First-Passage-Time Problem for a Wiener Process Annals of Applied
Probability, In press (2009)
\end{thebibliography}
\end{document}